\documentclass[aps,pra,superscriptaddress,twocolumn,longbibliography,floatfix,10pt]{revtex4-1}

\usepackage[normalem]{ulem}
\usepackage{amsfonts,amssymb,dsfont}
\usepackage[sumlimits,intlimits]{amsmath}
\usepackage{graphics}
\usepackage{graphicx}
\usepackage[dvipsnames]{xcolor}
\usepackage{color}
\usepackage{mathrsfs}
\usepackage{textcomp}
\usepackage{verbatim}
\usepackage{bm}
\usepackage{soul}
\usepackage{braket}
\usepackage{times}
\usepackage[T1]{fontenc} 
\usepackage{booktabs}
\usepackage{lipsum}
\usepackage{tabularx}

\newcommand{\T}[1]{\mathcal{T}\!{#1}}

\newcommand{\diver}[1]{\mathrm{div}\!\left({#1}\right)}

\usepackage{xr-hyper}
\usepackage[colorlinks=true,citecolor=blue]{hyperref}
\usepackage[capitalise,compress]{cleveref}

\newcommand*{\wh}{\widehat}

\newcommand{\Z}{\mathbb{Z}}
\newcommand{\R}{\mathbb{R}}
\newcommand{\C}{\mathbb{C}}
\newcommand{\N}{\mathbb{N}}
\newcommand{\dd}{\mathrm{d}}
\DeclareMathOperator{\erf}{erf}
\DeclareMathOperator{\erfc}{erfc}

\newcommand\scalemath[2]{\scalebox{#1}{\mbox{\ensuremath{\displaystyle #2}}}}

\begin{document}

	\title{Robust series linearization of nonlinear advection-diffusion equations}

	\author{T.~Forrest Kieffer}
	\thanks{Corresponding Author}
	\email{Thomas.Kieffer@jhuapl.edu}
	\affiliation{Johns Hopkins Applied Physics Laboratory, Laurel, Maryland 20723, USA}
	\author{Jakob Cupp}
	\affiliation{Johns Hopkins Applied Physics Laboratory, Laurel, Maryland 20723, USA}
	\author{John S.~Van Dyke}
	\affiliation{Johns Hopkins Applied Physics Laboratory, Laurel, Maryland 20723, USA}
	\author{Paraj Titum}
	\affiliation{Johns Hopkins Applied Physics Laboratory, Laurel, Maryland 20723, USA}
	\affiliation{William H. Miller III Department of Physics \& Astronomy, Johns Hopkins University, Baltimore, Maryland 21218, USA}
	\author{Michael L.~Wall}
	\affiliation{Johns Hopkins Applied Physics Laboratory, Laurel, Maryland 20723, USA}

	
\begin{abstract}	
We consider nonlinear partial differential equations (PDEs) for advection-diffusion processes which are augmented by an auxiliary parameter $\delta$ such that $\delta=0$ corresponds to linear advection-diffusion.  We derive potentially non-perturbative series expansions in $\delta$ that provide a process to obtain the solution of the nonlinear PDE through solving a hierarchical system of linear, forced PDEs with the forcing terms dependent on solutions at lower orders in the hierarchy. We rigorously detail our approach for a particular deformation that interpolates between linear advection-diffusion and the canonical Burgers' equation modeling nonlinear advection. In this case, we prove that the series has infinite radius of convergence for arbitrary integrable initial data, analyze the cases of a Dirac-delta initial condition (IC) (i.e., the fundamental solution) in an infinite domain and arbitrary IC in a periodic domain, and demonstrate the approach to turbulent behavior in a scenario with periodic forcing.  We then treat models of nonlinear diffusion involving the $p$-Laplacian operator, including generalizations of the Poisson equation in $1$ and $2$ dimensions, and the heat equation in $1+1$ dimensions.  We detail series expansions for two different deformations of these equations about their linear (ordinary Laplacian) counterparts, providing numerical evidence for the convergence of the series outside of a perturbative regime and demonstrating that the rate and radius of convergence are affected by choice of deformation.  Our results provide a rigorous foundation for using series expansion techniques to study nonlinear advection-diffusion PDEs, opening new pathways for analysis and potential applications for quantum-assisted computational fluid dynamics.
\end{abstract}
	
\maketitle
\section{Introduction}

Advection-diffusion equations, in which the concentration of a conserved quantity (e.g., particles or energy) is moved throughout a domain through transport (advection) and diffusion, occur ubiquitously in the physical sciences and contain many paradigmatic examples of parabolic partial differential equations (PDEs)~\cite{evans2022partial}.  In a general sense, we can parameterize advection-diffusion equations in $d$ dimensions as the following equation for a density $u : (0,\infty) \times \R^d\rightarrow \R$
\begin{align}
\label{eq:nlinPDE} \partial_t u  &= \diver{ \mathbf{j} (u,\nabla u) } \, ,
\end{align}
together with initial data $u(0,x) = g(x)$. We can interpret $\mathbf{j}$ as a flux of the conserved quantity, and several well-known instances of advection-diffusion follow from this interpretation. For example, if we consider the total flux to be the sum of an advective flux $\mathbf{j}_a=\mathbf{v} u$ with velocity $\mathbf{v}$ and a diffusive flux given by Fick's law $\mathbf{j}_d=-\nu \nabla u$ with $\nu$ the diffusivity (taken to be scalar for simplicity), we find the \emph{linear advection-diffusion equation}
\begin{align}
\label{eq:LADE} \partial_tu&=\mathbf{v}\cdot \nabla u-\nu \Delta u\, ,
\end{align}
where $\Delta = \sum_j \partial_j^2$ denotes the Laplace operator on $\R^d$. The linearity of this equation as a PDE is inherited from the total flux $\mathbf{j}=\mathbf{j}_a+\mathbf{j}_d$ being a linear function of $u$ and its gradient, but nonlinear generalizations are possible.  Specializing to $1$ dimension and considering the case in which the advective velocity is proportional to the concentration, ${j}_a=\frac{1}{2}u^2$, we obtain the canonical \emph{Burgers' equation} for nonlinear advection
\begin{align}
\label{eq:BE} \partial_tu&=u\nabla u-\nu \Delta u\, .
\end{align}
In a similar vein, if we consider a case in which the diffusivity depends on the concentration gradient to some power, we obtain the $p$-Laplacian evolution equation
\begin{align}\label{eq:pLaplacianEvolutionEqn}
\partial_t u = \Delta_p (u)
\end{align}
in which \begin{align}\label{eq:pLaplacian}
\Delta_p (u) = \diver{ |\nabla u|^{p-2} \nabla u } 
\end{align}
is the $p$-Laplacian operator, defined for $p\in\left(1,\infty\right)$, and $| \cdot |$ denotes the standard Euclidean norm of a vector in $\R^d$. When $p = 2$, the $p$-Laplacian is just the regular Laplace operator and \eqref{eq:pLaplacianEvolutionEqn} reduces to the heat equation for linear diffusion. 

Nonlinear advection occurs even for simple, Newtonian, incompressible fluids through the celebrated Navier-Stokes (NS) equations.  Burgers' equation Eq.~\eqref{eq:BE} arises as a simplification of NS through ignoring the pressure term, which additionally enables the equation to be consistently studied in one spatial dimension, in contrast to full NS.  Nonlinear diffusion also arises in fluid dynamics through large-eddy simulation (LES) models~\cite{guermond2004mathematical}, in which equations of fluid mechanics (e.g., NS) are filtered through a linear operator $u\to\bar{u}$ which commutes with constant-coefficient differential operators.  Qualitatively, the action of the filtering operator is to separate out the large-distance and small-distance scales, and so should produce a set of equations for the behavior at large scales with a source term feeding in the effects of smaller scales, perhaps phenomenologically.  Applying the filtering operator to the NS equations produces a system of equations that looks again like NS in $\bar{u}$, but with an additional diffusion term represented by the divergence of the the so-called subgrid-scale tensor $T\left(\nabla u\right)=\overline{\nabla u\otimes \nabla u}-\overline{\nabla u}\otimes\overline{\nabla u}$.  Ideally, the subgrid-scale tensor should be expressible as a function of $\bar{u}$ alone (i.e., make no reference to $u$) in order to simplify the problem, but such an \emph{exact closure} does not actually reduce the number of degrees of freedom, a situation called the \emph{closure paradox}~\cite{guermond2004mathematical}.  Hence, practical LES approaches look for accurate models for $T\left(\overline{\nabla u}\right)$ that are not expected to be exact but retain nice properties (e.g., that the filtered equations are provably well-posed).  The LES model of Smagorinsky~\cite{smagorinsky1963general} amounts to taking 
\begin{align}
T_{\varepsilon}\left(\nabla u\right)&=-\varepsilon^2 \left|D\right| D\, ,\\
D&=\left[\nabla u+\left(\nabla u\right)^T\right]/2\, ,
\end{align}
with $\varepsilon>0$ chosen to reproduce the physics of the underlying microscopic equations, such as the $k^{-5/3}$ scaling of energy-versus-wavenumber $k$ in NS.  The Smagorinsky model falls into the more general class of Ladyzenskaja models whose subgrid-scale diffusion tensor takes the form $T\left(\nabla u\right)=\beta\left(\left|\nabla u\right|^2\right)\nabla u$ where $\beta:[0,\infty)\to[0,\infty)$ is monotone increasing and $c_1x^{\mu}\le \beta(x) \le c_2 x^{\mu}$ for $\mu\ge 1/4$.  Ladyzenskaja models enjoy the property that they produce a well-posed system of equations that converge to weak solutions of NS in the limit $\varepsilon\to0$.  We also see that Ladyzenskaja models with $\beta$ a pure power law produce the $p$-Laplacian as a diffusion operator with $p\ge 2.5$ (Smagorinsky is the case $p=3$).  While these two examples highlight important cases where nonlinear advection-diffusion arises in fluid dynamics, they are by no means exhaustive; Burgers' equation appears in many other contexts to include condensed matter, nonlinear acoustics, and cosmology~\cite{frisch2002burgulence}, and the $p$-Laplacian appears in contexts of image processing, optimal transport, sandpile growth, and others~\cite{diening2013mini}.

Given the ubiquity of nonlinear advection-diffusion phenomena, there is considerable interest in solving nonlinear PDEs of the form Eq.~\eqref{eq:nlinPDE} and a range of methods have been developed.  Perhaps the simplest is direct numerical simulation (DNS), in which the solution is obtained in high spatial and temporal resolution without relying on additional models or uncontrolled approximations.  Finite difference or finite element methods on fine grids applied directly to the PDE are examples of DNS approaches.  However, as is well-known, phenomena such as turbulence introduce couplings between a very large range of scales, challenging the limits and the scalability of DNS applied to such problems.  Even in the case of effective nonlinear models which arise from treating the large-scale dynamics at high resolution and phenomenologically treating the fine-scale dynamics (e.g., the LES models described above), the scales involved can be stressing for DNS, and care must be taken that the numerical scheme captures the properties of the solution.  A relevant example arises for the finite-difference approximation to the $p$-Laplacian; the only known finite-difference approximation to the $p$-Laplacian which is monotone~\cite{del2022finite}, i.e. converges to the proper viscosity solution, involves an operator whose support $r$ scales non-trivially with $p$ and with the grid spacing $h$, requiring large stencils.  This can be contrasted with the ordinary Laplacian, where low-order stencils with $h\sim\mathcal{O}(r)$ are accurate and monotone.

Given the difficulties of DNS applied to PDEs with conventional computing, there has been significant interest in applying alternative computing paradigms such as quantum computing to the solution of PDEs, see Refs.~\cite{dalzell2023quantum,morales2024quantum, jaksch-et-al-2023-variational-quantum-algorithms-for-computational-fluid-dynamics,bharadwaj2020quantum,joseph2023quantum} for some recent reviews.  A significant consideration here is that quantum computations are expressed in terms of unitary operations and so inherently linear.  A linear PDE can be related to a linear system of equations through discretization as in DNS approaches, and then known quantum algorithms for linear systems such as the Harrow-Hassidim-Lloyd (HHL) algorithm~\cite{harrowquantum2009,krovi2023improved} can be applied.  Direct application of similar ideas to nonlinear PDEs requires multiple copies of a known fiducial state and incurs an exponential overhead~\cite{leyton2008quantum}, though this overhead can be reduced for certain types of polynomial nonlinearites together with variational means of representing quantum states~\cite{jaksch-et-al-2023-variational-quantum-algorithms-for-computational-fluid-dynamics}.  Hence, significant effort has been devoted to means of linearizing the PDE in a means amenable to a quantum solution, with Carleman~\cite{liu2021efficient} and Koopman-von Neumann~\cite{PhysRevResearch.2.043102} linearizations being prominent approaches.

In this work, we take an approach in which the nonlinear PDE is converted into an infinite hierarchy of linear PDEs by introducing a deformation of the equation controlled by a parameter $\delta$ such that $\delta=0$ is a linear equation and the nonlinear equation of interest is obtained at some value $\delta_{\star}$, expanding the solution of our equation as a Taylor series in $\delta$ about $\delta = 0$, and then matching coefficients in the series expansion through the PDE order-by-order.  As shown below, the resulting PDEs are all \emph{linear}, and form a hierarchy in which functionals of the solutions from lower order appear as forcing terms in the equations for higher-order terms.  This enables the collection of tools developed for linear PDEs to be brought to bear, including analytic techniques, local yet monotone finite difference schemes, and more straightforward adaptation for quantum computing.  

The idea of using series expansions applied to nonlinear PDEs is not new; Bender \emph{et al.} considered series solutions for a particular deformation of Burgers' equation in Ref.~\cite{bender1991new} and the related homotopy analysis method~\cite{liao2009notes} has a similar philosophy of using a homotopy to interpolate between a linear PDE and a nonlinear one.  Our work expands significantly beyond this body of research, demonstrating that the series approach does not have to be interpreted perturbatively through the explicit identification of a deformation for the viscous Burgers' equation proven to have infinite radius of convergence, and exploring the convergence behavior of nonlinear models of diffusion involving the $p$-Laplacian through nonlinear deformations centered around the Laplace operator.  These results provide a rigorous foundation for future studies and applications of series methods across the myriad application areas of nonlinear advection-diffusion equations.

Our work is organized as follows: in Sec.~\ref{sec:Burgers} we consider series methods applied to deformations of Burgers' equation, which features linear diffusion but nonlinear advection.  After presenting a methodology to obtain the hierarchy for a general homotopy connecting Burgers' to linear advection-diffusion, we study a particular deformation in detail, and show that the associated series has infinite radius of convergence given mild constraints on the IC.  We analyze the case of delta function IC to discuss practical matters of series convergence, and present an optimization we call \emph{refeeding} -- wherein a truncated series expansion is summed and used as a new IC for the series expansion technique -- that improves the convergence behavior.  We use these methodologies in a simulation of Burgers' equation with forcing on a periodic domain, and see hallmark characteristics of Burgers' turbulence~\cite{jeng1969forced,bec2007burgers}.  In Sec.~\ref{sec:pLap} we apply similar techniques to nonlinear diffusive models involving the $p$-Laplacian, studying both an evolution equation that generalizes the heat equation as well as a Dirichlet problem that generalizes Poisson's equation.  We again present a general method to derive a hierarchy of linear PDEs from a parameterized homotopy before studying two different homotopies -- which we call the ordinary and the dual-- in detail.  Using exactly known solutions for the $p$-Laplacian Dirichlet problem in 1D and 2D and the evolution equation in 1D, we study the convergence of these two series expansions for a range of $p$ and show validity outside of a naive perturbative range $p=2+\delta$, $\left|\delta\right|<1$.  Finally, in Sec.~\ref{sec:DaO}, we provide a discussion and outlook.  Some technical points of our exposition are presented as Appendices.

\subsubsection*{Notation} For $m \geq 0$ and $p \in [1,\infty]$, $W^{m,p} (\Omega)$ denotes the Sobolev space consisting of functions whose (distributional) derivatives up to order $m \geq 0$ are $p^{\mathrm{th}}$-power integrable over the (possibly unbounded) domain $\Omega \subseteq \R^d$. When $p = \infty$, $W^{m,\infty} (\Omega)$ denotes the Sobolev space consisting of functions whose (distributional) derivatives up to order $m \geq 0$ are essentially bounded on $\Omega \subseteq \R^d$. The space $W^{0,p} (\Omega)$ coincides with the standard Lebesgue space $L^p (\Omega) = \{ f : \Omega \rightarrow \R ~ : ~ \int_{\Omega} |f|^p < \infty \}$. Similarly, $L^{\infty} (\Omega)$ denotes the space of essentially bounded functions on $\Omega$. Most frequently, $|\cdot|$ will denote the standard Euclidean norm of a vector in $\R^d$. Occasionally we use $|\cdot|_{\C}$ to indicate the modulus of a complex number. 

\section{Nonlinear advection: deformations of Burgers' equation}\label{sec:Burgers}

We begin with a statement of the most general form of nonlinear advection we consider, referred to as the \emph{homotopy Burgers' equation} (HPE). On a possibly unbounded spatial interval $I \subset \R$, the HPE refers to the semilinear Cauchy problem
\begin{align}\label{eq:HomotopyBurgers}
\left\lbrace \begin{array}{l}
\partial_t u + \partial_x h(u,\delta) - \nu \partial^2_x u = f \\[5pt]
u (0,x) = g(x) ,
\end{array} \right. 
\end{align}
where $u : (0,\infty) \times I \times [0,1] \rightarrow \R$ is the unknown field considered as a function of $(t,x,\delta)$, $\nu > 0$ is the viscosity, $h : \R \times [0,1] \rightarrow \R$ is a smooth function that satisfies $h(u,0) = v u$, with $v \in \R$ a constant, and $h(u , 1) = u^2 / 2$, $f : (0,\infty) \times I \rightarrow \R$ is an external force, and $g : I \rightarrow \R$ is the given IC. If $I$ is bounded, then $u$ may have Dirichlet or periodic boundary conditions (BCs), while if $I = \R$, then we demand $u (t,x;\delta) \rightarrow 0$ as $x \rightarrow \pm \infty$ for any $(t, \delta) \in (0,\infty) \times [0,1]$ (and likewise for the initial data $g$). 

The conditions imposed on $h$ ensure that, as $\delta$ varies, Eq.~\eqref{eq:HomotopyBurgers} defines a family of nonlinear advection-diffusion equations which "deforms" a linear advection-diffusion equation ($\delta = 0$) into the classical Burgers' equation, Eq.~\eqref{eq:BE}, ($\delta = 1$). There are infinitely many such deformations; as an example, Ref.~\cite{bender1991new} considers $h(u,\delta) = u^{1 + \delta} / (1 + \delta)$ with $v=1$. The special case $h(u,\delta) = (1-\delta) v u + \delta u^2 / 2$, which we will call the \emph{linear homotopy}, will be considered in detail momentarily. It is possible (and perhaps desired) to also have $h$ depend explicitly on the independent variables $(t,x)$. This does not appreciably change the following exposition and so we ignore this  generalization.

Let $u (t,x;\delta)$ denote the solution to \eqref{eq:HomotopyBurgers}. Suppose $\delta \mapsto u (t,x;\delta)$ is analytic at $\delta = 0$ for almost every $(t,x) \in (0,\infty) \times I$. Then, there exists a sequence of functions $u_n : (0,\infty) \times I \rightarrow \R$ such that, for almost every $(t,x)$, $\delta \mapsto u (t,x;\delta)$ may be expanded as a convergent series with some radius of convergence $r = r (t,x) > 0$:
\begin{align}\label{eq:deltaExpansion}
u (t,x ; \delta) = \sum_{n \geq 0} \delta^n u_n (t,x) .
\end{align}
The functions $u_n$ are then computed via
\begin{align}\label{eq:un_def}
u_n (t,x) := \frac{1}{n!} \partial_{\delta}^n u (t,x ; \delta) \Big|_{\delta = 0} .
\end{align}
By expanding $h$ in a Taylor series about the point $(u_0 , 0)$, plugging this and \eqref{eq:deltaExpansion} into \eqref{eq:HomotopyBurgers}, and matching powers of $\delta$, we may derive a hierarchy of PDE each coefficient function $u_n$ must satisfy. The first few orders are derived here, with the general case treated in Appendix \ref{app:GeneralBE}.

We will suppress the dependence of $u$ on $x$ and $t$ to keep the notation uncluttered and write $u(\delta)$ for $u(t,x;\delta)$. First observe that $h(u (\delta) , \delta) |_{\delta = 0} = h(u_0 , 0) = v u_0$, which follows from the homotopy property of $h$. The chain rule gives that
\begin{align*}
\left. \frac{\dd}{\dd \delta} h(u (\delta) , \delta) \right|_{\delta = 0} &= h^{(1,0)}_0 u_1 + h^{(0,1)}_0 , \\[5pt]
\left. \frac{\dd^2}{\dd \delta^2} h(u (\delta) , \delta) \right|_{\delta = 0} &= h^{(2,0)}_0 u_1^2 + 2 h^{(1,1)}_0 u_1 \\
& \hspace{5mm} + h^{(0,2)}_0 + 2 h^{(1,0)}_0 u_2 ,
\end{align*}
where $h^{(n,m)}_0 := (\partial_1^n \partial_2^m h)(u_0,0)$ for integers $n, m \geq 0$. (Here, $\partial_1$ and $\partial_2$ indicate the partial derivatives with respect to the first and second arguments of the function, respectively.) Therefore, the first few equations in the hierarchy of PDE read
\begin{widetext}
\begin{align}\label{eq:BurgerHomotopyHierarchy}
\left\lbrace \begin{array}{l}
\partial_t u_0 + v \partial_x u_0 - \nu \partial_x^2 u_0 = f \\[2pt]
\partial_t u_1 + \partial_x \left( h^{(1,0)}_0 u_1 \right)  - \nu \partial_x^2 u_1 = - \partial_x h^{(0,1)}_0 \\[2pt]
\partial_t u_2 + \partial_x \left( h^{(1,0)}_0 u_2 \right)  - \nu \partial_x^2 u_2 = - \dfrac{1}{2} \partial_x \left( h^{(2,0)}_0 u_1^2 + 2 h^{(1,1)}_0 u_1 + h^{(0,2)}_0 \right) \\
\hspace{3cm} \vdots
\end{array}  \right. 
\end{align} 
\end{widetext}
These equations are supplemented with the ICs $u_0 (0,x) = g(x)$ and $u_n (0,x) = 0$ for $n \geq 1$. Any BCs that come along with \eqref{eq:HomotopyBurgers} may, in general, be handled in an analogous way as the ICs. 

An essential point is that each equation in \eqref{eq:BurgerHomotopyHierarchy} is \textit{linear} for the current order but with coefficients and forcing determined by $h$ and the lower orders. This pattern continues for each order, as shown explicitly in Appendix \ref{app:GeneralBE}. As examples, the choice $h(u,\delta) = u^{1 + \delta} / (1 + \delta)$ yields the first-order equation
\begin{align*}
\partial_t u_1 + \partial_x u_1  - \nu \partial_x^2 u_1  = - \ln{(u_0)} \partial_x u_0 ,
\end{align*}
in accordance with Ref.~\cite{bender1991new}, and the linear homotopy $h(u,\delta) = (1-\delta) v u + \delta u^2 / 2$ has first-order equation
\begin{align*}
\partial_t u_1 + v \partial_x u_1  - \nu \partial_x^2 u_1 = (v - u_0) \partial_x u_0 . 
\end{align*}
Clearly, the complexity of the advection and forcing term is dependent on the choice of homotopy.

As with any nonlinear PDE, well-posedness should not be taken for granted, and conditions should be imposed on $h$ to ensure that, for each $\delta \in (0,1)$, \eqref{eq:HomotopyBurgers} always possesses a unique time-global solution whenever $g$ and $f$ are appropriately regular. To appreciate this point, we note the classical result that the nonlinear heat equation $\partial_t u - \nu \partial_x^2 u = u^{\delta}$ equipped with an arbitrary non-negative, smooth, compactly supported initial data, does not possess a time-global, non-negative and integrable solution whenever $\delta \in (1,3)$ \cite{fujita1966blowing}. The remainder of this section will only be concerned with the simple choice of linear homotopy for $h$, for which time-global well-posedness is readily addressed through an explicit representation of the solution (c.f.~Eq.~\eqref{eq:BurgerLinearHomotopy_soln2}). Hence, it is outside the scope of this work to determine conditions on $h$ which ensure time-global well-posedness of Eq.~\eqref{eq:HomotopyBurgers}. We refer the interested reader to \cite{GUIDOLIN2022126361,BIANCHINI2025128761} for global existence and blow-up results for equations of the form \eqref{eq:HomotopyBurgers}, and \cite{Friedman1988,Bandle1994,BANDLE19983} for more information on well-posedness of semilinear diffusion equations in general.

Let $\mathcal{S}_N$ denote the partial sum of the series expansion \eqref{eq:deltaExpansion} up order $N$, i.e.
\begin{align}\label{eq:BurgersPartialSum}
\mathcal{S}_N ( t,x;\delta ) = \sum_{n=0}^N \delta^n u_n (t,x) ,
\end{align}
where each $u_n$ is obtained by solving each equation in \eqref{eq:BurgerHomotopyHierarchy}. Assume the solution of \eqref{eq:HomotopyBurgers} is analytic in $\delta$ at $\delta = 0$ for almost every $(t,x) \in (0,\infty) \times I$, and, furthermore, that the expansion \eqref{eq:deltaExpansion} converges at $\delta = 1$. Then, $\mathcal{S}_N ( t,x;1 )$ converges (point-wise almost everywhere in $(t,x)$) to the solution of Burgers' equation as $N \rightarrow \infty$. A key question of interest is then: what is the rate-of-convergence of $\mathcal{S}_N ( t,x;1 )$ as $N \rightarrow \infty$ as a function of $(t,x)$ and other problem parameters?

The remainder of this section is devoted to tackling this and closely related questions for the choice of linear homotopy $h$ under different BCs, initial data, and forcing. Section \ref{sec:LDSE} formally introduces the linear homotopy Burgers' equation and discusses its key properties. In particular, the Cole-Hopf transform~\cite{cole1951quasi,hopf1950partial} is used to show that $\C \ni \delta \mapsto u(t,x;\delta) \in \C$ is a holomorphic function for every $t > 0$ and \textit{almost every} $x \in \R$. One important practical implication of this result is that the radius of convergence of the Taylor series expansion about $\delta = 0$ \eqref{eq:deltaExpansion} is infinite for every $t > 0$ and almost every $x \in \R$, allowing us to take $\delta = 1$ in \eqref{eq:deltaExpansion}. Section \ref{sec:BurgersFS} will then look at the explicit solution of \eqref{eq:HomotopyBurgers} evolved from a Dirac delta-function IC. Here, we will study the rate-of-convergence of the series \eqref{eq:deltaExpansion} as a function of the variables $t$, $x$, and $v$ using the exact Taylor expansion coefficient functions $u_n (t,x)$. 

Section \ref{sec:PBCs} turns to \eqref{eq:HomotopyBurgers} on a compact spatial domain under periodic BCs.  Here, we will use spectral methods to numerically solve each equation in \eqref{eq:BurgerHomotopyHierarchy} up to a desired order, say $N$, in order to approximate the solution of Burgers' through the partial sum $\mathcal{S}_N ( t,x;1 )$. The solver will be validated using a closed-form solution to Burgers' corresponding to a cosine-squared IC and zero forcing. Then, we will introduce the important concept of refeeding, in which the truncated series expansion $\mathcal{S}_N ( t_f,x;1 )$ is obtained by integrating \eqref{eq:BurgerHomotopyHierarchy} over a small time interval $[0 , t_f]$ and reused as an IC to the first equation in \eqref{eq:BurgerHomotopyHierarchy}, upon which the hierarchy is evolved over another short time interval and the process is repeated. With refeeding, excellent convergence properties of the series expansion will be demonstrated. This section concludes with a study of Burgers' turbulence through the series expansion. In particular, we show that the turbulent energy-versus-wavenumber scaling characteristics of the forced Burgers' equation~\cite{jeng1969forced,bec2007burgers} is well captured using the series expansion technique. The main summarizing conclusion from this entire analysis is that solutions of Burgers' equation can be accurately approximated by solving an appropriate hierarchy of \textit{linear} advection-diffusion equations.

\subsection{The linear deformation, series expansion, and proof of analyticity}\label{sec:LDSE}

For $\delta \in \C$, consider the \textit{linear homotopy Burgers' equation}:
\begin{align}\label{eq:BurgerLinearHomotopy}
\left\lbrace \begin{array}{l}
\partial_t u + (1-\delta) v \partial_x u + \delta u \partial_x u - \nu \partial^2_x u = f \\[5pt]
u (0,x) = g (x) ,
\end{array} \right.
\end{align}
where, again, $v \in \R$ is an arbitrary constant dictating the velocity of the linear advection term. The choice of $v$ will impact the convergence rate of the series expansion and, using dimensional analysis, one may identify two "natural" choices for $v$. Let $\ell_0$, $u_0$, and $t_0 = \ell_0 / u_0$ be nominal length, velocity, and time scales, respectively. Define the dimensionless field $\tilde{u} (s,y)$ that is a function of dimensionless variables $(s,y) \in [0 , \infty) \times \R$ via $\tilde{u} (s,y) = u (t_0 s, \ell_0 y) / u_0$. The equation for $\tilde{u}$ then reads
\begin{align*}
\left\lbrace \begin{array}{l}
\partial_s \tilde{u} + (1-\delta) \dfrac{v}{u_0} \partial_y \tilde{u} + \delta \tilde{u} \partial_y \tilde{u} - \dfrac{1}{R_e} \partial^2_y \tilde{u} = \tilde{f} \\[5pt]
\tilde{u} (0,y) = \tilde{g} (y) .
\end{array} \right.
\end{align*}
where $\tilde{f} (s,y) = f(t_0 s, \ell_0 y)/u_0$, $\tilde{g} (y) = g(\ell_0 y) / u_0$, and $R_e = u_0 \ell_0 / \nu$ is the \textit{Reynolds number}. The two "natural" choices for $v$ are then $v = u_0$ {or} $v = \nu / \ell_0=u_0/R_e$. The latter choice will be shown to yield improved convergence properties of the expansion \eqref{eq:deltaExpansion} relative to the former choice. In any case, for the remainder of this manuscript, the linear homotopy Burgers' equation will refer to the dimensionless equation:
\begin{align}\label{eq:BurgerLinearHomotopy_dimensionless}
\left\lbrace \begin{array}{l}
\partial_t u + (1-\delta) v \partial_x u + \delta u \partial_x u - \dfrac{1}{R_e} \partial^2_x u = f \\[5pt]
u (0,x) = g (x) .
\end{array} \right.
\end{align}

The main problem we are concerned with here is demonstrating the analyticity of the solution $u$ to \eqref{eq:BurgerLinearHomotopy_dimensionless} as a function of the parameter $\delta \in \C$ in order to rigorously justify the use of the series expansion \eqref{eq:deltaExpansion}. To tackle this problem we specialize to the case of $I = \R$, $g \in L^1 (\R)$ real-valued, and $f = 0$, and argue that $\C \ni \delta \mapsto u (t,x;\delta) \in \C$ is holomorphic at $\delta = 0$ through the Cole-Hopf transform. Then, we consider the radius of convergence $0 < r(t,x) \leq \infty$ of the Taylor expansion of $\delta \mapsto u (t,x;\delta)$ about $\delta = 0$ and argue that $r(t,x) = \infty$ for every $t > 0$ and \textit{almost every} $x \in \R$. The important implication of this result is that the Taylor series expansion of $\delta \mapsto u(t,x;\delta)$ about $\delta = 0$, i.e. \eqref{eq:deltaExpansion}, is convergent for every $\delta \in \C$ and, in particular, at $\delta = 1$. Equivalently, the functions $\{ u_n \}_{n \geq 0}$ that satisfy \eqref{eq:BurgerHomotopyHierarchy} with the linear homotopy (see Eq.~\eqref{eq:LinearHomotopyHierarchy} for an explicit expression) yield a convergent series $\sum_{n \geq 0} \delta^n u_n (t,x)$ for every $\delta \in \C$ and which solves the original $\delta$-dependent nonlinear PDE. This result provides rigorous justification for using the series expansion to solve Burgers' equation (i.e., taking $\delta = 1$), at least in the special case when $f = 0$.

The Cole-Hopf transform refers to a function $w : (0,\infty) \times \R \times \C \rightarrow \C$ that is related to the solution $u$ of Eq.~\eqref{eq:BurgerLinearHomotopy_dimensionless} via
\begin{align}\label{eq:ColeHopf}
u = - \frac{2}{\delta R_e} \partial_x \ln{w} .
\end{align}
Note that $w$ is only uniquely defined up to an arbitrary time-dependent multiplicative factor. It may be readily verified that $w$ satisfies the following linear PDE:
\begin{align}\label{eq:BurgerLinearHomotopy_ColeHopf}
\left\lbrace \begin{array}{l}
\partial_t w + (1-\delta) v \partial_x w - \dfrac{1}{R_e} \partial_x^2 w = -\dfrac{\delta R_e}{2} (\T{f}) w \\[5pt]
w (0,x) = \exp{ \left( -\frac{\delta R_e}{2} (\T{g}) (x) \right) } ,
\end{array} \right.
\end{align}
where $\mathcal{T} : L^1 (\R) \rightarrow C(\R) \cap L^{\infty} (\R)$ is the linear integral operator defined by
\begin{align}\label{eq:Toperator}
(\T{g}) (x) := \int_{-\infty}^x g (y) \dd y .
\end{align}

Let's introduce the linear differential operator
\begin{align}\label{eq:Loperator}
\mathcal{L}_{\delta} := (\delta - 1) v \partial_x + \frac{1}{R_e} \partial_x^2 
\end{align}
and let $e^{t \mathcal{L}_{\delta}}$ denote the associated evolution semigroup for $t > 0$. The action of $e^{t \mathcal{L}_{\delta}}$ on a function $f \in L^{\infty} (\R)$ is the convolution of a Gaussian integral kernel against $f$:
\begin{align}\label{eq:Ldelta_Semigroup}
(e^{t \mathcal{L}_{\delta}} f) (x) = \sqrt{ \frac{R_e}{4 \pi t} } \int_{\R} e^{- \frac{R_e (x-y-(1-\delta)vt)^2}{4 t}} f(y) \dd y . 
\end{align}
In this notation, the unique solution to \eqref{eq:BurgerLinearHomotopy_ColeHopf} when $f = 0$ is provided by the formula
\begin{align}\label{eq:BurgerLinearHomotopy_ColeHopf_soln}
w(t,x;\delta) = (e^{t \mathcal{L}_{\delta}} e^{ -\frac{\delta R_e}{2} \T{g} }) (x) .
\end{align}
Note that Eq.~\eqref{eq:BurgerLinearHomotopy_ColeHopf_soln} makes sense for any $\delta = \delta_1 + i \delta_2 \in \C$ because, for $g \in L^1 (\R)$, $e^{ -\frac{\delta R_e}{2} \T{g} } \in L^{\infty} (\R)$, which implies $w (t, \cdot ; \delta) \in L^{\infty} (\R)$.

With \eqref{eq:BurgerLinearHomotopy_ColeHopf_soln}, we may transform back into a solution for $u$ using \eqref{eq:ColeHopf}. The resulting general solution of \eqref{eq:BurgerLinearHomotopy_dimensionless} with $f = 0$ reads
\begin{align}\label{eq:BurgerLinearHomotopy_soln}
u (t,x ; \delta) = - \frac{2}{\delta R_e} \partial_x \ln{ (e^{t \mathcal{L}_{\delta}} e^{ -\frac{\delta R_e}{2} \T{g} }) (x) } .
\end{align}
Expression \eqref{eq:BurgerLinearHomotopy_soln} is valid for any $\delta \in \R \backslash \{0\}$ since the argument in the logarithm is positive for any $\delta \in \R \backslash \{0\}$. In particular, \eqref{eq:BurgerLinearHomotopy_soln} gives the correct solution to Burgers' when evaluated at $\delta = 1$ and to linear advection-diffusion when $\delta\to 0 $ following an application of L'H$\hat{\mathrm{o}}$pital's rule. These observations justify carrying out the differentiation in \eqref{eq:BurgerLinearHomotopy_soln} and using $[\partial_x , \mathcal{L}_{\delta}]=0$ to write
\begin{align}\label{eq:BurgerLinearHomotopy_soln2}
u (t,x ; \delta) = \frac{(e^{t \mathcal{L}_{\delta}} g e^{ -\frac{\delta R_e}{2} \T{g} })(x)}{(e^{t \mathcal{L}_{\delta}} e^{ -\frac{\delta R_e}{2} \T{g} })(x)} .
\end{align}

The numerator and denominator in \eqref{eq:BurgerLinearHomotopy_soln2} are separately seen to be holomorphic functions of $\delta \in \C$ because they each involve the composition, product, and integrals involving holomorphic functions of $\delta$. Hence, the ratio in \eqref{eq:BurgerLinearHomotopy_soln2} defines a meromorphic function of $\delta$. In particular, the function $\delta \mapsto u(t,x;\delta)$ is holomorphic at $\delta = 0$ since the denominator is $1$ there. The radius of convergence of the series expansion of $\delta \mapsto u(t,x;\delta)$ about $\delta = 0$ is determined by the distance to the nearest non-removable singularity. Singularities, if any, will arise as zeros of the denominator in \eqref{eq:BurgerLinearHomotopy_soln2}. Hence, a desirable assertion would be that $(e^{t \mathcal{L}_{\delta}} e^{ -\frac{\delta R_e}{2} \T{g} })(x) \neq 0$ for any $(t,x,\delta) \in (0,\infty) \times \R \times \C$ whenever $g \in L^1 (\R)$. 

It is straightforward to verify that the denominator in \eqref{eq:BurgerLinearHomotopy_soln2} never vanishes for $\delta \in \R$ because $(e^{t \mathcal{L}_{\delta}} e^{ -\frac{\delta R_e}{2} \T{g} })(x) > 0$ in this case. Hence, the challenge is ruling out the existence of a $\delta \in \C$ with non-zero imaginary part for which there is a $(t,x) \in (0,\infty) \times \R$ such that $(e^{t \mathcal{L}_{\delta}} e^{ -\frac{\delta R_e}{2} \T{g} })(x) = 0$. Now, note that $\R \times \C \ni (x,\delta) \mapsto e^{ -\frac{\delta R_e}{2} (\T{g})(x) } \in \C$ is never zero because $\T{g} \in L^{\infty} (\R)$. A preliminary question is then whether $|(e^{\frac{t}{R_e} \partial_x^2}  e^{ -\frac{\delta R_e}{2} \T{g} }) (x)|_{\C} > 0$ for all $t > 0$ and $x \in \R$, where $|\cdot|_{\C}$ indicates the complex modulus. 

We have not succeeded in demonstrating this strict inequality \textit{for all} $x \in \R$ given our stated assumptions on $g$. In fact, it is plausible that evolution of $e^{ -\frac{\delta R_e}{2} \T{g} }$ under the heat kernel $e^{\frac{t}{R_e} \partial_x^2}$ produces a zero at some $(t,x) \in (0,\infty) \in \R$, even though $|e^{ -\frac{\delta R_e}{2} \T{g} (x) }|_{\C} > 0$ \footnote{The more general claim that, given $f : \R \rightarrow \C$, $|f|_{\C} > 0$ implies $|e^{t \partial_x^2} f|_{\C} > 0$ is \textit{false}. A counterexample we know of is a superposition of the form $f (x) = ( c_1 e^{i \xi_1 x} + c_2 e^{i \xi_2 x} ) e^{- x^2 / 2}$ for appropriately chosen $c_1 , c_2 , \xi_1 , \xi_2 \in \R$. However, this counterexample is not of the form $e^{ -\frac{\delta R_e}{2} \T{g} (x) }$, and so doesn't immediately imply that the assertion $|(e^{\frac{t}{R_e} \partial_x^2}  e^{ -\frac{\delta R_e}{2} \T{g} }) (x)|_{\C} > 0$ is false.}. However, $\C \ni x \mapsto (e^{\frac{t}{R_e} \partial_x^2}  e^{ -\frac{\delta R_e}{2} \T{g} }) (x) \in \C$ is readily seen to be an entire function, and entire functions cannot vanish on a subset of $\C$ that contain a limit point as a consequence of the identity theorem for holomorphic functions. Consequently, the set of $x \in \R$ for which $(e^{t \mathcal{L}_{\delta}} e^{ -\frac{\delta R_e}{2} \T{g} })(x) = (e^{\frac{t}{R_e} \partial_x^2}  e^{ -\frac{\delta R_e}{2} \T{g} }) (x - (1-\delta)vt)$ vanishes must have zero measure. This, in turn, implies the radius of convergence of the Taylor expansion of $\delta \mapsto u(t,x;\delta)$ about $\delta = 0$ is infinite for every $t > 0$ and almost every $x \in \R$. 


The above outlined proof of analyticity demonstrates that the series expansion technique exemplified in Eq.~\eqref{eq:BurgerHomotopyHierarchy} is not necessarily perturbative; i.e. $|\delta|$ is not restricted to small values in order that the process converge. One caveat here is that we restricted ourselves to specific BCs and $f = 0$. Rigorous analysis in the case of $f \neq 0$ appears far more delicate. For instance, it is not even immediately clear what assumptions must be made on $f$ in order to ensure the solution of Eq.~\eqref{eq:BurgerLinearHomotopy_ColeHopf} is holomorphic at $\delta = 0$, let alone what additional restrictions must be placed on $f$ for the radius of convergence of the Taylor expansion of $\delta \mapsto u(t,x;\delta)$ about $\delta = 0$ to include $\delta = 1$. Despite the lack of rigorous results for $f \neq 0$, numerical evidence in the following subsections strongly suggest analyticity at $\delta = 0$ and convergence of \eqref{eq:deltaExpansion} at $\delta = 1$ in situations with $f \neq 0$. Questions regarding analyticity with respect to parameters for fundamental solutions of constant coefficient \textit{linear} PDE goes back to at least H\"{o}rmander \cite{HormanderALPDOII}, but the mathematical literature for similar questions regarding nonlinear PDE (or even non-constant coefficient linear PDE) appears sparse.


\subsubsection{Case study: The fundamental solution}
\label{sec:BurgersFS}

As an application of our techniques, consider the case of a Dirac-delta IC supplied to \eqref{eq:BurgerLinearHomotopy_dimensionless} with $f=0$. We refer to the corresponding time-evolved field as the \emph{fundamental solution}. From \eqref{eq:BurgerLinearHomotopy_soln}, this fundamental solution reads
\begin{align*}
\vspace{-0.5\baselineskip}
\scalemath{0.9}{
u (t,x ; \delta) = - \frac{2}{\delta R_e} \partial_x \ln{\int_{\R} e^{- \frac{R_e (x-y-(1-\delta)t)^2}{4 t} - \frac{\delta R_e}{2} \mathbf{1}_{(0,\infty)} (y)} \dd y} ,
}
\end{align*}
where $\mathbf{1}_A$ denotes the indicator function for the set $A \subseteq \R$. Performing the integration and simplifying results in
\begin{align}\label{eq:BurgerLinearHomotopy_deltaIC}
\scalemath{0.9}{
u (t,x ; \delta) = \frac{2 (e^{\delta R_e} - 1) e^{- \frac{(x - (1-\delta) v t)^2}{4 t / R_e}} }{ \delta \sqrt{\pi R_e t} \left( 2 + (e^{\delta R_e} - 1) \erfc{\left( \frac{x - (1-\delta) v t}{2 \sqrt{t / R_e}} \right)} \right) } .
}
\end{align}
Equation \eqref{eq:BurgerLinearHomotopy_deltaIC} reproduces the correct Burgers' solution for $\delta = 1$~\cite{evans2022partial}. Moreover, 
\begin{align}\label{eq:BurgerLinearHomotopy_deltaIC_u0}
u_0 (t,x) := \lim_{\delta \rightarrow 0} u (t,x ; \delta) = \sqrt{ \frac{ R_e }{ \pi t } } e^{- \frac{R_e (x-vt)^2}{4 t} } ,
\end{align}
which satisfies $\partial_t u_0 + v \partial_x u_0 - \frac{1}{R_e} \partial_x^2 u_0 = 0$. The first-order term $u_1$ in the Taylor expansion of \eqref{eq:BurgerLinearHomotopy_deltaIC} about $\delta = 0$, i.e. $u_1 (t,x) = \lim_{\delta \rightarrow 0} \partial_{\delta} u(t,x;\delta)$, is given by
\begin{align*}
\scalemath{0.9}{
u_1 (t,x)  = \sqrt{ \frac{ R_e^3 }{ 4 \pi t } } e^{- \frac{R_e (x-vt)^2}{4 t} } \left( \erf{\left( \frac{x-vt}{2 \sqrt{t/R_e}} \right)} - v(x - vt) \right)
}
\end{align*}
Together $u_0$ and $u_1$ satisfy $\partial_t u_1 + v \partial_x u_1 - \frac{1}{R_e} \partial_x^2 u_1 = (v - u_0) \partial_x u_0$.

\begin{figure}[h!]
  \begin{center}
\includegraphics[width=1.0\columnwidth]{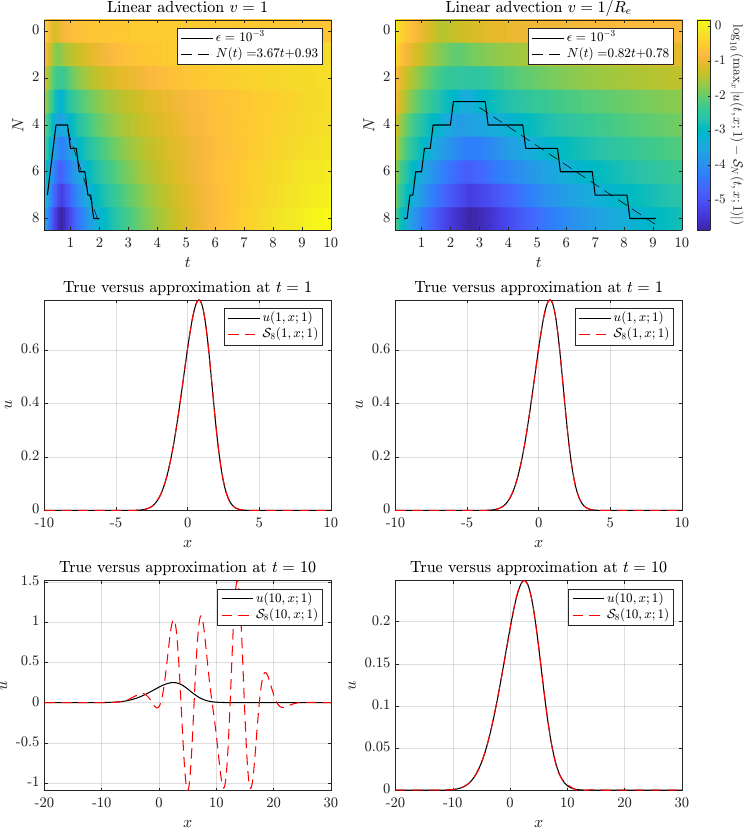}  
\caption{\label{fig:Burger_deltaIC}  \emph{Convergence of series expansion for the fundamental solution of Burgers' equation.} The exact solution of Burgers' corresponding to a Dirac-delta IC (i.e., \eqref{eq:BurgerLinearHomotopy_deltaIC} with $\delta = 1$) and $R_e = 2$, is compared to the exact partial sum $\mathcal{S}_N (t,x;1)$ for varying expansion order $N$, time $t$, and (normalized) linear advection speed $v$. Left and right columns correspond to $v = 1$ and $v = 1/R_e$, respectively. Panels in the top row show the maximum spatial error between $u (t,x;1)$ and $\mathcal{S}_N (t,x;1)$ for $0 \leq N \leq 8$ and $0 < t \leq 10$. The solid black lines overlaying these plots indicate the error threshold of $\epsilon = 10^{-3}$, while the dashed black lines indicate a least-squares linear fit to a portion of this data. Middle and bottom rows show the exact solution compared to $\mathcal{S}_8$ for $t = 1$ and $t = 10$, respectively.}
\end{center}
\end{figure}

With the aid of computer algebra software (e.g., Mathematica), it is possible to exactly compute $u_n$ to arbitrary order using \eqref{eq:BurgerLinearHomotopy_deltaIC} and compare the resulting truncated Taylor expansion $\mathcal{S}_N$ \eqref{eq:BurgersPartialSum} to the true solution of Burgers' equation (i.e., \eqref{eq:BurgerLinearHomotopy_deltaIC} evaluated at $\delta = 1$). Figure \ref{fig:Burger_deltaIC} captures such a comparison for $R_e = 2$ and varying expansion order, time, and (normalized) linear advection speed $v$. A modest Reynolds number was chosen here to avoid numerical issues when evaluating the exact Taylor expansion coefficient functions $u_n$. Higher Reynolds number will be considered in subsequent sections. 

Figure \ref{fig:Burger_deltaIC} illustrates that $\mathcal{S}_N$ is an excellent approximation of the true solution to Burgers' equation for $v = 1/R_e$ and $N \geq 4$. The choice $v = 1$ also has relatively low error for $t \lesssim 4$, but the error in $\mathcal{S}_N$ at fixed order $N$ begins to diverge as $t$ grows. Increasing the order $N$ of the expansion increases the time duration with which one can remain under a reasonable error threshold. The same phenomenon may be observed for the choice $v = 1/R_e$, although it is not as dramatic. By studying the exact expression \eqref{eq:BurgerLinearHomotopy_deltaIC}, we know that $\delta \mapsto u(t,x;\delta)$ is analytic at $\delta = 0$ with infinite radius of convergence. Thus, given an error threshold $\epsilon > 0$, it is always possible to find an $N = N(t,\epsilon)$ such that $\max_x | u(t,x;1) - \mathcal{S}_N (t,x;1) | < \epsilon$. However, the numerical evidence in Figure \ref{fig:Burger_deltaIC} strongly suggests that $N \rightarrow \infty$ as $t \rightarrow \infty$. Hence, it is not practical to use this series expansion technique if one aims to study long-term behavior because, in its current state, the technique requires one to solve many linear advection-diffusion equations \eqref{eq:BurgerHomotopyHierarchy} over long time intervals. Fortunately, this issue will be remedied in the following section using the concept of refeeding. 

\subsection{Burgers' turbulence through the series expansion}
\label{sec:PBCs}

This section will consider \eqref{eq:BurgerLinearHomotopy_dimensionless} on the spatial domain $[0,1] \subset \R$ with periodic BCs. In other words, we seek the solution $u : (0,\infty) \times [0,1] \times \C \rightarrow \C$ to \eqref{eq:BurgerLinearHomotopy_dimensionless} which satisfies $u(t,0;\delta) = u(t,1;\delta)$ for all $\delta \in \C$ and $t > 0$. In particular, the initial data $g$ must satisfy $g(0) = g(1)$. The Cole-Hopf transform \eqref{eq:ColeHopf} may again be employed to transform the periodic version of \eqref{eq:BurgerLinearHomotopy_dimensionless} into \eqref{eq:BurgerLinearHomotopy_ColeHopf}, but with the definition of $\mathcal{T}$ now modified to read
\begin{align}
(\T{g}) (x) := \int_0^x g (y) \dd y .
\end{align}
The aim is to leverage the (circular) Fourier transform to solve the periodic version of \eqref{eq:BurgerLinearHomotopy_ColeHopf}. However, the IC in \eqref{eq:BurgerLinearHomotopy_ColeHopf}, i.e. $\exp{ \left( -\frac{\delta R_e}{2} (\T{g}) (x) \right) }$, will not be periodic unless $(\T{g})(1) = \int_0^1 g = 0$. This condition can always be ensured by shifting $g \mapsto g - \int_0^1 g$ and so we consider this condition fulfilled in the following discussion.

The (periodic) Fourier transform of a function $h \in C^1 (\mathbb{S}^1)$ on the circle $\mathbb{S}^1 = \R / \Z$ is denoted by 
\begin{align}\label{eq:periodicFT}
\wh{h} (k) = \int_0^1 e^{- 2\pi i k x} h(x) \dd x .
\end{align}
Note $h \in C^1 (\mathbb{S}^1)$ ensures $\wh{h} \in \ell^1 (\Z)$. Fourier inversion in this case gives the classical Fourier series
\begin{align}\label{eq:FourierSeries}
h(x) = \sum_{k \in \Z} \wh{h} (k) e^{2 \pi i k x} . 
\end{align}
Recalling that $\wh{ \partial_x h } (k) = 2 \pi i k \wh{h} (k)$ for $h \in C^1 (\mathbb{S}^1)$, the Fourier transform of \eqref{eq:BurgerLinearHomotopy_ColeHopf} when $f = 0$ reads
\begin{align*}
\partial_t \wh{w} (k) + 2 \pi i k (1-\delta) v \wh{w} (k) + \frac{4 \pi^2 k^2}{R_e} \wh{w} (k) = 0 ,
\end{align*}
which is readily solved to find
\begin{align*}
\wh{w} (t,k) = e^{- (2 \pi i k (1-\delta) v + 4 \pi^2 k^2 / R_e) t} \wh{w} (0,k) .
\end{align*}
Inverting with \eqref{eq:FourierSeries} then yields
\begin{align}\label{eq:BurgerLinearHomotopy_periodic_ColeHopf_soln}
\scalemath{0.9}{
w(t,x) = \int_0^1 \left( \sum_{k \in \Z} e^{2\pi i k (x - (1-\delta) v t - y) - \frac{4 \pi^2 k^2 t}{R_e}} \right) w(0,y) \dd y .
}
\end{align}

Using the Poisson summation formula, the inner series in \eqref{eq:BurgerLinearHomotopy_periodic_ColeHopf_soln} defines a Gaussian distribution that has been wrapped to the circle, i.e.
\begin{align}\label{eq:WrappedNormal}
\scalemath{0.9}{
\Phi_{2 + \delta}^t (x) = \sqrt{ \frac{R_e}{4 \pi t} } \sum_{k \in \Z} \exp{ \left( - \frac{R_e (x - (1-\delta) v t - k)^2}{4 t} \right) } .
}
\end{align}
In terms of \eqref{eq:WrappedNormal}, we can write the solution to our periodic linear advection-diffusion equation \eqref{eq:BurgerLinearHomotopy_ColeHopf} as the convolution
\begin{align*}
w(t,x) = (\Phi^t_{\delta} \ast w(0,\cdot)) (x) := \int_0^1 \Phi_{2 + \delta}^t (x-y) w(0,y) \dd y .
\end{align*}
Returning to the Cole-Hopf transform \eqref{eq:ColeHopf}, the solution to the periodic linear homotopy Burgers' equation \eqref{eq:BurgerLinearHomotopy_dimensionless} reads
\begin{align*}
u(t,x;\delta) = - \frac{2}{\delta R_e} \partial_x \ln{ \left( \Phi_{2 + \delta}^t \ast e^{ -\frac{\delta R_e}{2} \T{g} } \right) (x) } \, ,
\end{align*}
or, equivalently,
\begin{align}\label{eq:BurgerLinearHomotopy_periodic_soln}
u(t,x;\delta) = \frac{ (\Phi_{2 + \delta}^t \ast g e^{ -\frac{\delta R_e}{2} \T{g} }) (x) }{ (\Phi_{2 + \delta}^t \ast e^{ -\frac{\delta R_e}{2} \T{g} } ) (x) }  .
\end{align}
Notice the similarity between \eqref{eq:BurgerLinearHomotopy_periodic_soln} and \eqref{eq:BurgerLinearHomotopy_soln2}.

As an application of the preceding analysis, and to develop a closed-form solution for benchmarking, consider the initial data $g(x) = \cos^2{(2 \pi x)} - 1/2$. To develop a closed-form solution for this IC, we rely on the Jacobi-Anger expansion:
\begin{align*}
e^{iz\sin{\theta}} = \sum_{k \in \Z} J_k (z) e^{ik\theta} ,
\end{align*}
valid for $\theta \in \C$, where $J_k : \C \rightarrow \C$ is the Bessel function of the first kind of order $k \in \Z$. Since $(\T{g})(x) = \sin{(4\pi x)} / 8\pi$, we may apply the Jacobi-Anger expansion to arrive at the identity
\begin{align}\label{eq:JacobiAngerExpnsion}
e^{ -\frac{\delta R_e}{2} (\T{g}) (x) } = \sum_{k \in \Z} I_k \left( \frac{\delta R_e}{16 \pi} \right) e^{ 2 \pi i k \left( 2 x + \frac{1}{4} \right) } ,
\end{align}
where $I_k : \C \rightarrow \C$ is the modified Bessel function of the first kind of order $k \in \Z$. Identity \eqref{eq:JacobiAngerExpnsion} together with formula \eqref{eq:BurgerLinearHomotopy_periodic_soln} yields
\begin{widetext}
\begin{align}\label{eq:BurgerLinearHomotopy_periodic_soln_cosinesquaredIC}
u(t,x;\delta) = \frac{16 \pi}{\delta R_e} \partial_x \ln{ \left( I_0 \left( \frac{\delta R_e}{16 \pi} \right) + 2 \sum_{k \geq 1} I_k \left( \frac{\delta R_e}{16 \pi} \right) \cos{ \left( 4 \pi k \left( x - (1-\delta) v t + \frac{1}{8} \right) \right) }  e^{- 16 \pi^2 k^2 t / R_e} \right) } .
\end{align}
\end{widetext}
As with the Dirac-delta IC from Section \ref{sec:BurgersFS}, one may use \eqref{eq:BurgerLinearHomotopy_periodic_soln_cosinesquaredIC} to compute, for example, $u_0 (t,x) = \lim_{\delta \rightarrow 0} u(t,x;\delta)$ and $u_1 (t,x) = \lim_{\delta \rightarrow 0} \partial_{\delta} u(t,x;\delta)$, and verify that these functions satisfy $\partial_t u_0 + v \partial_x u_0 - \frac{1}{R_e} \partial_x^2 u_0 = 0$ with $u_0 (0,x) = \cos^2{(2 \pi x)} - 1/2$ and $\partial_t u_1 + v \partial_x u_1 - \frac{1}{R_e} \partial_x^2 u_1 = (v - u_0) \partial_x u_0$ with $u_1 (0,x) = 0$.

In general, the equations for $u_n$, $n \geq 1$, satisfy (see Eq.~\eqref{eq:LinearHomotopyHierarchy})
\begin{align*}
\scalemath{0.9}{
\left\lbrace \begin{array}{l}
\partial_t \wh{u}_n (t,k) + \left( 2 \pi i k v + \frac{ 4 \pi^2 k^2 }{ R_e } t \right) \wh{u}_n (t,k) = - 2 \pi i k \wh{F}_n (t,k) \\[5pt]
\wh{u}_n (0,k)  = 0 ,
\end{array} \right.
}
\end{align*}
where 
\begin{align}
\scalemath{0.9}{
\wh{F}_n (t,k) := \int_0^1 e^{-2 \pi i k x} F_n (u_0 (t,x) , \cdots , u_{n-1} (t,x)) \dd x ,
}
\end{align}
and $F_n (u_0 (t,x) , \cdots , u_{n-1} (t,x))$ is given by \eqref{eq:InhomogeneityLinearHomotopy} in the appendix. Using Duhamel's formula and Fourier inversion, we arrive at a convenient representation formula for $u_n$:
\begin{align}\label{eq:BurgerLinearHomotopy_periodic_uk}
\scalemath{0.825}{
u_n (t,x) = - 2 \pi i \sum_{k \in \Z} k e^{2 \pi i k x} \int_0^t e^{- \left( 2 \pi i k v + \frac{ 4 \pi^2 k^2 }{R_e} \right) (t-\tau) } \wh{F}_n (\tau,k) \dd \tau .
}
\end{align}
Formula \eqref{eq:BurgerLinearHomotopy_periodic_uk} is the core of our fast Fourier transform (FFT)-based PDE solver for the hierarchy \eqref{eq:LinearHomotopyHierarchy} in the case of periodic BCs.

\begin{figure}[h!]
  \begin{center}
\includegraphics[width=1.0\columnwidth]{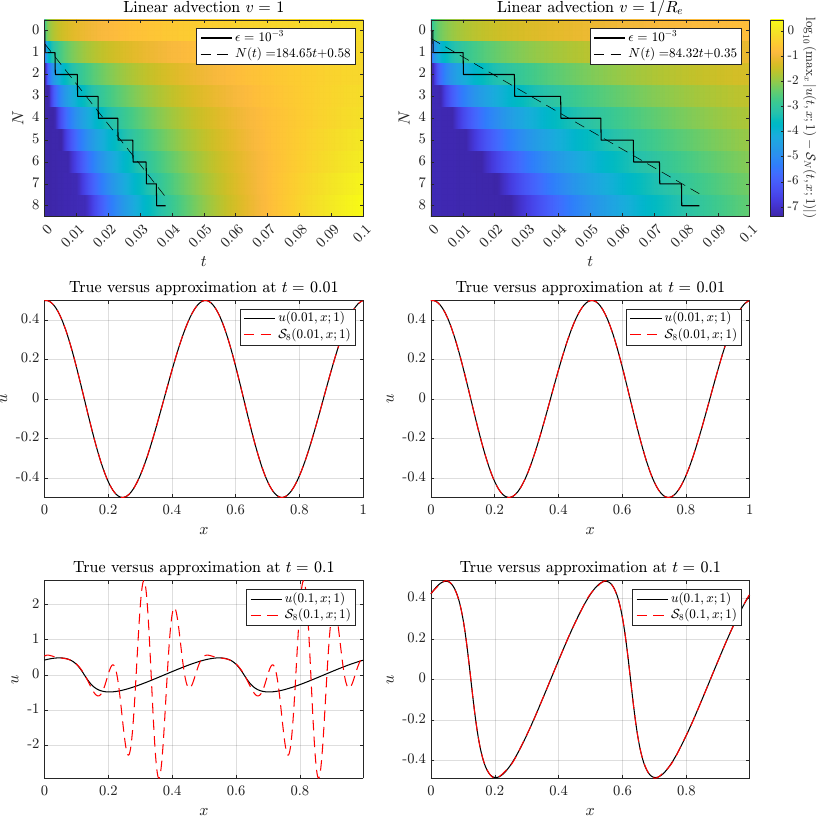}  
\caption{\label{fig:Burger_CosineSquaredIC}  \emph{Convergence of series expansion for the periodic Burgers' equation.} The exact solution of Burgers' corresponding to a cosine-squared IC (i.e., Eq.~\eqref{eq:BurgerLinearHomotopy_periodic_soln_cosinesquaredIC} with $\delta = 1$) and $R_e = 500$, is compared to the exact partial sum $\mathcal{S}_N (t,x;1)$ for varying expansion order $N$, time $t$, and (normalized) linear advection speed $v$. Left and right columns correspond to $v = 1$ and $v = 1/R_e$, respectively. Panels in the top row show the maximum spatial error between $u (t,x;1)$ and $\mathcal{S}_N (t,x;1)$ for $0 \leq N \leq 8$ and $0 < t \leq 0.1$. The solid black lines overlaying these plots indicate the error threshold of $\epsilon = 10^{-3}$, while the dashed black lines indicate a least-squares linear fit to a portion of this data. Middle and bottom rows show the exact solution compared to $\mathcal{S}_8$ for $t = 0.01$ and $t = 0.1$, respectively.}
\end{center}
\end{figure}

Figure \ref{fig:Burger_CosineSquaredIC} is similar to Figure \ref{fig:Burger_deltaIC}, but now using the exact solution corresponding to a cosine-squared IC, leveraging the FFT to numerically solve \eqref{eq:BurgerLinearHomotopy_periodic_uk}, and setting $R_e = 500$. Our grid spacing is set to $2^9 = 512$ (implying that the range of wavenumbers considered for the FFT spans $-2^8, -2^8 + 1 , \cdots , 2^8 - 1$) and the integration time step is $10^{-4}$. Again, the choice of (normalized) linear advection speed $v = 1/R_e$ is seen to be superior relative to $v = 1$, and both choices show diverging error as $t$ increases at fixed order. This divergence will now be remedied. 

\subsubsection{Refeeding}
\label{sec:refeeding}
Consider solving each equation in \eqref{eq:BurgerHomotopyHierarchy} over a time interval $[0,t_0]$, $t_0 > 0$, up to some desired order $N$ and constructing the partial sum $\mathcal{S}_N (t_0 , x ; \delta)$ \eqref{eq:BurgersPartialSum} at the time $t_0$. One may then begin to solve each equation in \eqref{eq:BurgerHomotopyHierarchy} over the next time interval $[t_0,t_1]$, $t_1 > t_0$, but now with initial data $u_0 (t_0,x) = \mathcal{S}_N (t_0 , x ; \delta)$ and $u_n (t_0 , x) = 0$, $n\ge 1$. This process is then repeated until the entire desired time interval is exhausted. We refer to this iterative process as refeeding.  

\begin{figure*}
  \begin{center}
\includegraphics[width=2\columnwidth]{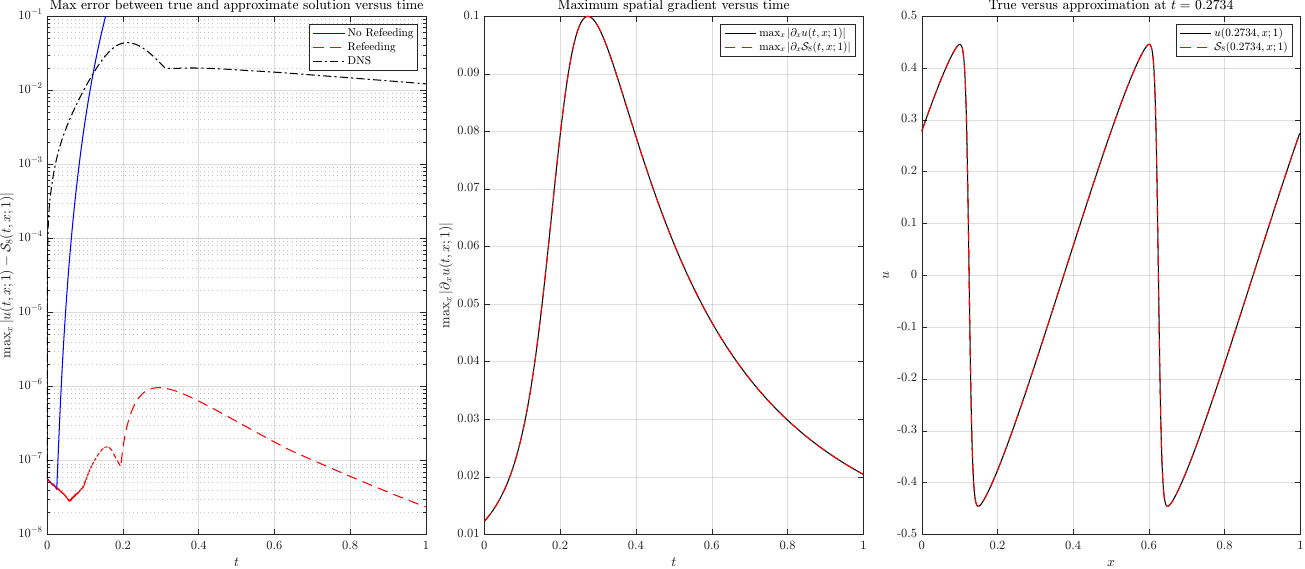}  
\caption{\label{fig:Burger_CosineSquaredIC_RefeedingComparison}  \emph{Convergence when refeeding is employed for periodic Burgers' equation.} The exact solution of Burgers' corresponding to a cosine-squared IC (i.e., Eq.~\eqref{eq:BurgerLinearHomotopy_periodic_soln_cosinesquaredIC} with $\delta = 1$) and $R_e = 500$, is compared to the approximate series solution with and without refeeding for expansion order $N = 8$ and (normalized) linear advection speed $v = 1/R_e$. Left panel shows the maximum spatial error between $u (t,x;1)$ and $\mathcal{S}_N (t,x;1)$ versus $t$ with and without refeeding alongside a standard DNS approach to solving Burgers' equation. Middle panel shows the maximum spatial gradient of both the exact solution and $\mathcal{S}_8$ versus time. Right panel shows the exact solution compared to $\mathcal{S}_8$ with refeeding at the time where the $\max_x |\partial_x u (t,x;1)|$ is largest.}
\end{center}
\end{figure*}

Figure \ref{fig:Burger_CosineSquaredIC_RefeedingComparison} shows the error in $\mathcal{S}_8$ versus time with and without refeeding alongside a standard DNS approach for Burgers' equation. The parameters used for this comparison are the same as those used to generate Figure \ref{fig:Burger_CosineSquaredIC}, with the exception that the final time considered is increased from $t = 0.1$ to $t = 1$ in order to highlight how the error remains small for long times when refeeding is employed. The refeeding time step is taken to be the same as the integration time step: $10^{-4}$. The DNS approach used for comparison is based on centered finite difference for the spatial derivatives and a forward Euler method for the time stepping. Figure \ref{fig:Burger_CosineSquaredIC_RefeedingComparison} shows how the series approach together with refeeding produces solutions which are orders-of-magnitude more accurate compared to a more standard numerical method. This is primarily due to the linear nature of each PDE in the hierarchy \eqref{eq:BurgerHomotopyHierarchy} which enabled a direct application of spectral methods, yielding both a fast and accurate solver. This highlights a key advantage of the series approach to the numerical simulation of nonlinear advection-diffusion equations. 

Figure \ref{fig:Burger_CosineSquaredIC_RefeedingComparison} also shows the maximum spatial gradient of both the true solution and $\mathcal{S}_8$ versus time, as well as the true and approximate solutions themselves at the time $t \simeq 0.2734$, which is when the maximum spatial gradient is maximized. Note that the place of greatest error, i.e. difference between $\mathcal{S}_8$ and the true solution, correlates with the time of the maximum spatial gradient.  It is evident from Figure \ref{fig:Burger_CosineSquaredIC_RefeedingComparison} that refeeding performs exceedingly well when numerically solving Burgers' equation using the series expansion $\mathcal{S}_N$ for the choice of linear homotopy. The success of refeeding is naturally related to continuous dependence of solutions to Burgers' on initial data in, e.g., $L^1$.

\begin{figure}
  \begin{center}
\includegraphics[width=0.9\columnwidth]{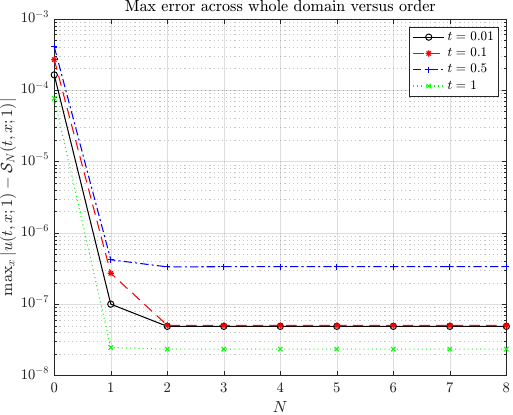}  
\caption{\label{fig:Burger_CosineSquaredIC_ErrVersusOrder} \emph{Maximum spatial error versus expansion order for periodic Burgers' equation.} The maximum spatial error between $u (t,x;1)$ and $\mathcal{S}_N (t,x;1)$ versus $N$ with refeeding employed at the time slices $t = 0.01$, $0.1$, $0.5$, and $1$.}
\end{center}
\end{figure}

Figure \ref{fig:Burger_CosineSquaredIC_ErrVersusOrder} shows the error versus expansion order $N$ for several time slices when refeeding is employed. The scenario parameters used here are the same as those used to generate the refeeding case from Figure \ref{fig:Burger_CosineSquaredIC_RefeedingComparison}. It is apparent from Figure \ref{fig:Burger_CosineSquaredIC_ErrVersusOrder} that the error rapidly converges, reaching a plateau once $N \geq 2$. This behavior may be understood by analyzing the first few solutions $u_n$ from \eqref{eq:LinearHomotopyHierarchy}. Following an application of Duhamel, each $u_0 (t) = e^{t \mathcal{L}_0} g$ and, for $n \geq 1$,
\begin{align}\label{eq:BurgerLinearHomotopyHeirarchy_Integral}
u_n (t) = \int_0^t e^{(t-s) \mathcal{L}_0} F_n (u_0 (s) , \cdots , u_{n-1} (s)) \dd s .
\end{align}
From the midpoint rule applied to each order in \eqref{eq:BurgerLinearHomotopyHeirarchy_Integral} and the polynomial nature of the $n^{\mathrm{th}}$-inhomogeneity $F_n$ \eqref{eq:InhomogeneityLinearHomotopy} in the variables $u_0, \cdots , u_n$, we conclude $u_n (t) \simeq \mathcal{O} (t^n)$. Therefore, if the integration time step for refeeding is small (e.g., $\sim 10^{-4}$), then there will be rapid convergence of $\mathcal{S}_N$ as a function of $N$. 

\subsubsection{Turbulence}

As an example application of the series approach to solving Burgers' equation using the linear homotopy, we consider the periodic version of \eqref{eq:BurgerLinearHomotopy_dimensionless} on the spatial domain $[0,1] \subset \R$ with forcing given by
\begin{align}\label{eq:BurgersPeriodicForcing}
f (x) = \sum_{k = k_{\mathrm{min}}}^{k_{\mathrm{max}}} A_k \sin{(2 \pi k x + \phi_k)} ,
\end{align}
where $A_k \in [-1,1]$ and $\phi_k \in [0,2\pi)$ are a set of amplitudes and phases, respectively, selected from uniform distributions over their respective domains. An exact, closed-form solution of \eqref{eq:BurgerLinearHomotopy_dimensionless} when $\delta = 1$ corresponding to this choice of forcing appears to be unknown. (See \cite{okamura1983steady} for a discussion of steady-state solutions when $k_{\mathrm{min}} = k_{\mathrm{\max}} = 1$. Indeed, after applying Cole-Hopf, Hill's equation results and some progress can be made by applying Floquet theory.) Our goal here is to use the series approach to numerically approximate the resulting turbulent solution of Burgers' through $\mathcal{S}_N (t,x;1)$. We will demonstrate that $\mathcal{S}_N (t,x;1)$ possesses the correct turbulent energy-versus-wavenumber scaling in the steady-state. 

\begin{figure*}
  \begin{center}
\includegraphics[width=2\columnwidth]{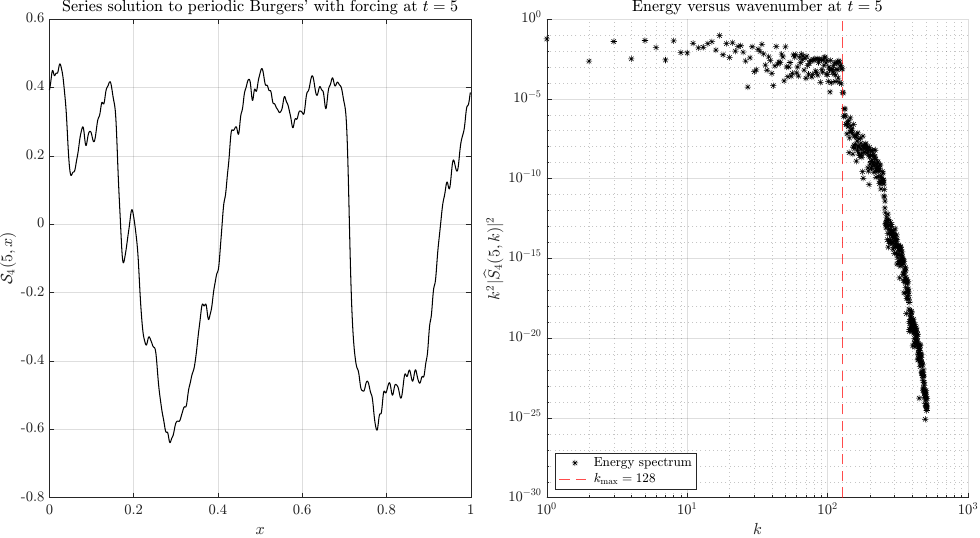}  
\caption{\label{fig:TurbulentSolnBurgers} \emph{Turbulent steady-state solution of the periodic Burgers' equation.} Left panel shows the turbulent steady-state solution of the periodic Burgers' equation with $R_e = 500$ and sinusoidal forcing \eqref{eq:BurgersPeriodicForcing} with $k_{\mathrm{min}} = 1$, $k_{\mathrm{max}} = 128$, and randomly sampled amplitudes and phases. The right panel shows $k^2 E(t,k)$ versus $k$ to illustrate the numerically-obtained turbulent steady-state exhibits the correct $E(t,k) \sim k^{-2}$ scaling in the steady-state.}
\end{center}
\end{figure*}

Consider the time-dependent \textit{kinetic energy functional} associated with the (periodic) Burgers' equation:
\begin{align}
\mathcal{E} (t) = \frac{1}{2} \int_0^1 |u(t,x)|^2 \dd x . 
\end{align}
Following an application of Parseval's identity, one may write $\mathcal{E} (t)$ as
\begin{align}
\mathcal{E} (t) = \frac{1}{2} \sum_{k \in \Z} |\wh{u} (t,k)|_{\C}^2 .
\end{align}
Since $u$ is real-valued, the instantaneous kinetic energy $E$ associated with wavenumber $k \geq 0$ is taken to be
\begin{align}
E (t,k) = |\wh{u} (t,k)|_{\C}^2 .
\end{align}
It has long been known \cite{jeng1969forced} that a turbulent steady-state solution of Burgers' satisfies the scaling $E(t,k) \sim k^{-2}$ in the inertial range (this is analogous to Kolmogorov's famous $k^{-5/3}$ scaling law for the Naiver-Stokes equations). Here, the lower and upper limit of the inertial range is determined by $k_{\mathrm{min}}$ and $k_{\mathrm{max}}$ in \eqref{eq:BurgersPeriodicForcing}, respectively. 

Figure \ref{fig:TurbulentSolnBurgers} shows a turbulent steady-state solution of Burgers' equation alongside a plot of $k^2 E(t,k)$ versus $k$ corresponding to the choices $R_e = 500$, $k_{\mathrm{min}} = 1$, and $k_{\mathrm{max}} = 128$. We've produced this turbulent steady-state by numerically solving the hierarchy up to $4^{\mathrm{th}}$ order (with refeeding) and then constructing the partial sum $\mathcal{S}_4 (5,x;1)$. The solution was evolved to $t = 5$ to ensure it has reached a steady-state. As expected, Figure \ref{fig:TurbulentSolnBurgers} shows that the scaled energy spectrum $k^2 E(t,k)$ is roughly flat for $1 \leq k \leq 128 = k_{\mathrm{max}}$, after which viscous dissipation determines the shape of the energy spectrum. 

\section{Nonlinear diffusion: deformations involving the $p$-Laplacian}
\label{sec:pLap}

We now turn our attention to nonlinear diffusion processes centered around the $p$-Laplacian operator, Eq.~\eqref{eq:pLaplacian}. Similar to how the homotopy Burgers' equation, Eq.~\eqref{eq:HomotopyBurgers}, interpolates between linear and nonlinear advection-diffusion through a choice of deformation $H(u,\delta)$, we define the \emph{homotopy $p$-Laplacian evolution equation}
\begin{align}\label{eq:HomotopyEvolution_pLaplacian}
\left\lbrace \begin{array}{l}
\partial_t u = \diver{|\nabla u|^{h (\delta)} \nabla u} \\[5pt]
u (0,x) = g(x) ,
\end{array} \right.
\end{align}
where $u : (0,\infty) \times \Omega \times [0,\delta_p] \rightarrow \R$ is the unknown field defined on a (possibly unbounded) spatial domain $\Omega \subseteq \R^d$ considered as a function of $(t,x,\delta)$, $g : \R^d \rightarrow \R$ is the IC, and $h : [0,\delta_p] \rightarrow \R$ is a smooth homotopy that satisfies $h (0) = 0$ and $h (\delta_p) = p-2$. With this choice, the right-hand side of Eq.~\eqref{eq:HomotopyEvolution_pLaplacian} becomes the ordinary $2$-Laplacian acting on $u$ at $\delta=0$ and becomes the $p$-Laplacian at $\delta = \delta_p$. 

As we did in Eqs.~\eqref{eq:deltaExpansion}-\eqref{eq:un_def}, suppose the map $\delta \mapsto u (t,x;\delta)$, where $u$ is the solution to \eqref{eq:HomotopyEvolution_pLaplacian}, is analytic at $\delta = 0$ for almost every $(t,x) \in (0,\infty) \times \Omega$. By considering the Taylor expansion of $\delta \mapsto \partial_t u (\delta) - \diver{|\nabla u (\delta)|^{h (\delta)} \nabla u (\delta)}$ about $\delta = 0$ and matching powers of $\delta$ in Eq.~\eqref{eq:HomotopyEvolution_pLaplacian}, we can find a hierarchy of PDEs each $u_n = \frac{1}{n!} \lim_{\delta \rightarrow 0} \partial_{\delta}^n u (\delta)$ must satisfy. The general case is treated in Appendix \ref{app:seriespLap}; here we only list the first three orders:
\begin{widetext}
\begin{align}\label{eq:pLaplacian_Hierarchy_generalHomotopy}
\left\lbrace \begin{array}{l}
(\partial_t - \Delta) u_0 = 0 \, ,\\[5pt]
(\partial_t - \Delta) u_1 = h' (0) \diver{ \nabla u_0 \ln \left|\nabla u_0 \right| }\, ,\\[5pt]
(\partial_t - \Delta) u_2 = \diver{ h' (0) \nabla u_1 \ln\left|\nabla u_0\right|  + \dfrac{1}{2} \nabla u_0 \left( h' (0) \ln^2\left|\nabla u_0\right| + 2 h' (0) \frac{ \langle \nabla u_0,\nabla u_1\rangle }{\left|\nabla u_0\right|^2} + h'' (0) \ln\left|\nabla u_0\right| \right) } \, ,\\[5pt]
\hspace{5cm} \vdots
\end{array} \right.
\end{align}
\end{widetext}
with corresponding ICs $u_0 (0,x) = g(x)$ and $u_n (0,x) = 0$ for $n \geq 1$. Similar to what we saw in \eqref{eq:BurgerHomotopyHierarchy}, each equation in \eqref{eq:pLaplacian_Hierarchy_generalHomotopy} is a linear, inhomogeneous heat equation in the unknown function $u_n$:
\begin{align}\label{eq:un_diff}
(\partial_t - \Delta) u_n = \diver{F_n (\nabla u_0, \cdots , \nabla u_{n-1})} ,
\end{align}
with a forcing term $F_n$ dependent only on lower order functions (see Appendix \ref{app:seriespLap} for an explicit expression). In terms of these inhomogeneities we can provide explicit solutions to the equations. For $n = 0$, we simply have
\begin{align}\label{eq:u0Int}
u_0 (t) = e^{t \Delta} g ,
\end{align}
where $e^{t \Delta}$ denotes the heat semigroup,
\begin{align}\label{eq:HeatSemiGroup}
(e^{t \Delta} g) (x) := \int_{\R^d} H_2 (t,x-y) g(y) \dd y \, ,
\end{align}
$H_2 (t,x)$ is the Gaussian heat kernel
\begin{align}\label{eq:HeatKernel}
H_2(t,x) = \frac{e^{- |x|^2 / 4t}}{(4 \pi t)^{d/2}} \, ,
\end{align}
and $g$ is the IC from \eqref{eq:HomotopyEvolution_pLaplacian}. For $n \geq 1$, the solution is given using Duhamel's principle as (c.f.~Eq.~\eqref{eq:BurgerLinearHomotopyHeirarchy_Integral})
\begin{align}\label{eq:un_int}
\scalemath{0.9}{
u_n (t) = \int_0^t e^{(t-s) \Delta} \diver{F_n (\nabla u_0 (s), \cdots , \nabla u_{n-1} (s))} \dd s .
}
\end{align}
For example, the solution to the equation for $u_1$ in \eqref{eq:pLaplacian_Hierarchy_generalHomotopy} reads
\begin{align}
u_1 (t) = \int_0^t e^{(t-s) \Delta} \diver{ \ln{|\nabla e^{s \Delta} g|} \nabla e^{s \Delta} g } \dd s .
\end{align}

In addition to the evolution equation Eq.~\eqref{eq:pLaplacianEvolutionEqn}, we will also consider the static $p$-Laplacian Dirichlet problem on a bounded domain $\Omega \subset \R^d$:
\begin{align}\label{eq:pLaplacianDirichlet}
\left\lbrace \begin{array}{ll}
\Delta_p (u) = f & \text{in } \Omega \\[5pt]
u = g  & \text{on } \partial \Omega  ,
\end{array} \right.
\end{align}
where $f, g : \Omega \rightarrow \R$ are given a priori. Eq.~\eqref{eq:pLaplacianDirichlet} generalizes the Dirichlet problem for Poisson's equation, which it reduces to when $p=2$. Analogously to the homotopy $p$-Laplacian evolution equation, we can define the homotopy $p$-Laplacian Dirichlet problem as
\begin{align}
\label{eq:HomotopyDirichlet_pLaplacian}
\left\lbrace \begin{array}{ll}
\diver{|\nabla u|^{h (\delta)} \nabla u} = f & \text{in } \Omega \, ,\\[5pt]
u = g  & \text{on } \partial  \Omega  ,
\end{array} \right.
\end{align}
where $h$ was defined in Eq.~\eqref{eq:HomotopyEvolution_pLaplacian}. We emphasize that the restriction to Dirichlet BCs is not a necessity and our ideas may be readily adapted to other types of BCs. An expansion very similar to Eq.~\eqref{eq:pLaplacian_Hierarchy_generalHomotopy} is available for the Dirichlet problem Eq.~\eqref{eq:HomotopyDirichlet_pLaplacian}.


PDE involving the $p$-Laplace operator must generally be interpreted in a weak sense because of the degenerate and singular nature of the operator for $p > 2$ and $1 < p < 2$, respectively. For example, the Dirichlet problem \eqref{eq:pLaplacianDirichlet} may be interpreted as the Euler-Lagrange equation for the functional $I_p : C(\Omega) \cap W^{1,p} (\Omega) \rightarrow \R$ given by
\begin{align*}
I_p (u) = \int_{\Omega} \left( f(x) u(x) - \frac{1}{p} |\nabla u (x)|^p \right) \dd x 
\end{align*}
with boundary data $g \in C (\overline{\Omega})$, $u - g \in W^{1,p}_0 (\Omega)$, and $f \in L^{\infty} (\Omega)$. A unique minimizer of this variational problem exists, and is actually $C^{1 , \alpha}_{\mathrm{loc}} (\Omega)$ with $\alpha \in (0,1)$, but, unless $p = 2$, the solution is generally not $C^2$ at points where the gradient vanishes \cite{lindqvist2019notes}, implying \eqref{eq:pLaplacianDirichlet} cannot be interpreted pointwise. Likewise, there are a number of results concerning the $p$-Laplacian evolution equation. Perhaps the most relevant to our exposition is the existence of a unique $u \in C((0,\infty) , L^1 (\R^d))$ whenever $p > 2$, $g \in L^1 (\R^d)$, and $u$ is stipulated to satisfy a growth condition at infinity \cite{dibenedetto1989cauchy}. Again, $C^2$-regularity of solutions is absent. 

An issue related to well-posedness is how solutions to PDEs involving the $p$-Laplacian depend on $p$. Under what conditions do we have (real) analyticity with respect to $p$? Moreover, what the is nature of the radius of convergence for the expansion of the solutions about $p = 2$? Unlike our results regarding the analyticity of solutions to homotopy Burgers' equation \eqref{eq:HomotopyBurgers} as a function of homotopy parameter $\delta$, this analogous question for PDEs involving the $p$-Laplacian appears far more challenging.  Here, it is important to clarify what is meant by analyticity here. We do not necessarily mean that the solution $u(x;p)$ to, say, the $p$-Laplacian Dirichlet problem \eqref{eq:pLaplacianDirichlet} satisfies
\begin{align*}
\lim_{N \rightarrow \infty} \sum_{n = 0}^N (p-2)^n u_n (x) = u (x;p)
\end{align*}
pointwise for almost every $x \in \Omega$, as what ended up being the case for the similar series expansion for the homotopy Burgers' equation \eqref{eq:HomotopyBurgers}. Instead, it is more plausible that the partial sums $\sum_{n = 0}^N (p-2)^n u_n (x)$ converge to the minimizer of $I_p (u)$ in, say, $L^p (\Omega)$. We will abuse the language, however, and not distinguish between these two situations when using the terminology "analytic with respect to $p$".

For the Dirichlet problem, this question of $p$-analyticity is more naturally phrased as whether minimizers of the functional $I_p$ are (real) analytic in $p$. In $d = 1$, this is readily answered in the positive by applying well-established results on the parameter dependence of solutions to nonlinear ordinary differential equations (ODEs) \cite{Coddington1956theory}. In fact, the Dirichlet problem \eqref{eq:pLaplacianDirichlet} with $d = 1$ is simply a boundary value problem for a nonlinear ODE:
\begin{align}\label{eq:pLaplacianDirichlet_1D}
\left\lbrace \begin{array}{ll}
(|u'|^{p-2} u')' = f & \text{in } (a,b) \subset \R \\[5pt]
u(a), u(b)  & \text{given}.
\end{array} \right.
\end{align}
Let $J_p : \R \rightarrow \R$ denote the bijective function $J_p (x) = |x|^{p-2} x$. Then we may write the ODE \eqref{eq:pLaplacianDirichlet_1D} as
\begin{align}
u' (x) = J_p^{-1} \left( |u'(a)|^{p-2} u' (a) + \int_a^x f (y) \dd y \right) . 
\end{align}
Upon another integration, we have the solution to the ODE in \eqref{eq:pLaplacianDirichlet_1D} in terms of the data $u(a)$ and $u'(a)$. This can of course be transformed into a solution of the boundary value problem \eqref{eq:pLaplacianDirichlet_1D} by solving an equation for $u'(a)$ in terms of $u(b)$. In any case, one can deduce from the explicit representation of $u$ in terms of $f$ and the boundary data that this solution is analytic for $p \in (1,\infty)$. This argument may be extended to radial solutions of \eqref{eq:pLaplacianDirichlet} when the domain is a ball in $\R^d$ and the source $f$ possesses radial symmetry. As a final remark, it is worth mentioning that when $d = 2$ and $f = 0$, the hodograph transform turns the Dirichlet problem \eqref{eq:pLaplacianDirichlet} into a linear PDE with non-constant coefficients (see \cite{aronsson1988p} for details). We suspect analyticity may be deduced from this observation, but it is not pursued further here. 

In general, only much weaker results are known rigorously. For example, the author in \cite{LINDQVIST198793} considers the stability of the Dirichlet problem \eqref{eq:pLaplacianDirichlet} with respect to $p$ by showing that $\lim_{p \rightarrow q} \| \nabla u (\cdot; p) - \nabla u (\cdot; q) \|_{L^p (\Omega)}^p = 0$ (under some suitable assumptions). Similar results available for the $p$-Laplacian evolution equation (see, e.g., \cite{KinnunenParviainen2010}). The challenge here is due, in part, to the lack of an explicit representation of solutions whose $p$-dependence may be directly analyzed. Indeed, unlike how Cole-Hopf transforms Burgers' equation into a linear PDE, there is no (known) transform which turns the $p$-Laplacian evolution equation or Dirichlet problem into a linear PDE. (As already mentioned, the exception is hodograph transform applied to the $p$-harmonic equation in $d = 2$.) This prevents us from obtaining representations of solutions with explicit $p$-dependence for general initial data in the evolution case and general sources, domains, and boundary data in the Dirichlet case. However, there are a few special cases for which closed-form solutions with explicit $p$-dependence are known, and these solutions are readily seen to be analytic with respect to $p$. We will employ some of these closed-form solutions in the sequel in order to gain insight into the properties and convergence behavior of the series expansions.

The remainder of this section will be concerned with numerically demonstrating the convergence of the series expansion and studying its convergence properties for both $p$-Laplacian evolution and Dirichlet problems. Section \ref{sec:pLaplacianOrdinaryDual} begins with a discussion of the two primary choices of homotopy $h$ used to generate series expansions. Section \ref{sec:pLaplaceDirichlet} focuses on the $p$-Laplace Dirichlet problem in $d = 1$ and $d = 2$. Special attention is given to a known closed-form solution on the ball and our ability to approximate this function with the partial sum $\sum_{n = 0}^N (p-2)^n u_n$, where $u_n$ is obtained by numerically solving its corresponding linear Poisson equation in the hierarchy. Section \ref{sec:pLaplaceEvolution} carries out a similar analysis for the $p$-Laplacian evolution equation. 

\subsection{The ordinary and dual deformations of the $p$-Laplacian and their series expansions}\label{sec:pLaplacianOrdinaryDual}

In this section some explicit forms for the homotopy function $h (\delta)$ are considered. Perhaps the simplest is the homotopy $h (\delta) = \delta$, in which case $\delta_p = p-2$. In other words, $p-2$ is treated as a "small parameter" in problem \eqref{eq:HomotopyEvolution_pLaplacian} and the solution is expanded in a power series about $p = 2$. The first three orders in the series expansion of Eq.~\eqref{eq:HomotopyEvolution_pLaplacian} with this deformation follow immediately from \eqref{eq:pLaplacian_Hierarchy_generalHomotopy} with $h'(0) = 1$ and $h''(0) = 0$. Similarly, the same homotopy applied to the static problem Eq.~\eqref{eq:HomotopyDirichlet_pLaplacian} yields a hierarchy of linear Dirichlet problems with sources at each order given by similar expressions as the inhomogeneities in \eqref{eq:pLaplacian_Hierarchy_generalHomotopy}. In this case, we supply $u_0$ with the same BC as in \eqref{eq:pLaplacianDirichlet}, while all other $u_n$ are zero on the boundary. We will refer to this homotopy and the resulting series expansions as \emph{ordinary}.



Naively, we can expect the ordinary series to converge when $|\delta|<1$, which in turn implies convergence for $1<p<3$.  As noted in the introduction, the $p$-Laplacian as an operator only makes sense (as written) for $p\in [1,\infty)$, and so we may further expect that $|\delta|<1$ may be required for convergence on the grounds that the operator makes sense. However, we may be interested in values of $p\ge 3$ for which convergence is less clear \emph{a priori}. A similar ordinary expansion could be performed for \eqref{eq:HomotopyEvolution_pLaplacian} with $h(\delta)$ there replaced by, say, $h(\delta) + 1$ and the domain of convergence overlap with our desired value of $p$. However, we lose the benefit of the resulting equations in the hierarchy being linear. 

These considerations motivate an alternate homotopy, which we will refer to as the \emph{dual}. We will arrive at the homotopy by considering the notion of the H\"{o}lder conjugate $p'$ of $1 \leq p \leq \infty$, defined by
\begin{align}\label{eq:HolderConjugate}
\frac{1}{p} + \frac{1}{p'} = 1 ,
\end{align}
or, equivalently, $p = \frac{p'}{p' - 1}$ and $p' = \frac{p}{p-1}$. Note that the H\"{o}lder conjugate operation maps $[1,\infty]$ to itself, and $p \geq 2$ if and only if $p' \leq 2$, with $p = p' = 2$ being the fixed point. In particular, if $p \geq 3$, then $p' \leq 3/2$. One significant aspect of $p'$ comes from the fact that the dual of the Lebesgue space $L^p (\Omega)$ is $L^{p'} (\Omega)$ for $1 < p < \infty$. $L^2$ is self-dual, a key feature of the Hilbert space aspect of $L^2$.

Suppose $p' = 2 + \delta$, where $|\delta| < 1$ is a small parameter. Then, $p-2 = (2-p')/(p'-1) = - \delta / (1+\delta)$. Notice how $3/2 < p < 3$ if and only if $3/2 < p' < 3$. This implies that both $p'$ and $p$ can be written as $2 + \delta$ for $|\delta| < 1$ when $p,p' \in (3/2 , 3)$. However, if $p \geq 3$, then only $p' = 2 + \delta$ for $|\delta| < 1$. For example, if $p = 4$, then $p' = 4/3 = 2 - 2/3$ and $\delta = - 2/3$.  It is this property of being able to explore large values of $p$ while keeping $\delta$ ``small" (i.e., keeping the $p$-Laplace operator close to the linear Laplacian) that motivates the choice of dual homotopy
\begin{align}\label{eq:dualhomotopy}
h (\delta) = - \frac{\delta}{1+\delta} ,
\end{align}
with $\delta_p = (2-p)/(p-1)$, by analogy with the H\"{o}lder dual to the ordinary homotopy. The first few orders in the dual series are obtained from Eq.~\eqref{eq:pLaplacian_Hierarchy_generalHomotopy} with $h'(0) = -1$ and $h'' (0) = 2$. A similar dual expansion can also be obtained for the Dirichlet problem.


\subsection{Convergence analysis - Dirichlet problem}\label{sec:pLaplaceDirichlet}

\begin{figure}
  \begin{center}
\includegraphics[width=1\columnwidth]{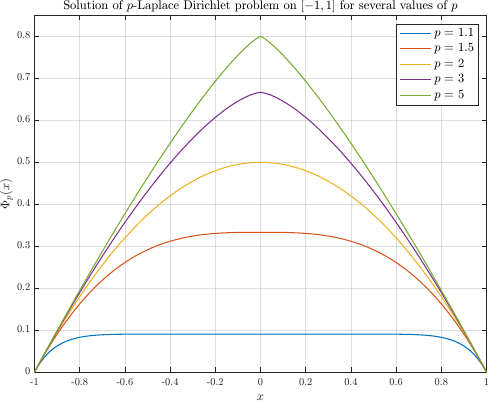}  
\caption{\label{fig:fundamentalSolution}  \emph{Exact solution of the $p$-Laplacian Dirichlet problem Eq.~\eqref{eq:fundamentalSolution}.} The exact solution to the $p$-Laplacian Dirichlet problem, Eq.~\eqref{eq:fundamentalSolution}, is shown for a range of $p$. Note that the solution is only smooth for $p \leq 2$; values of $p > 2$ feature a cusp at the origin.}
\end{center}
\end{figure}

Consider the Dirichlet problem \eqref{eq:pLaplacianDirichlet} on the unit ball $B_1 (0) \subset \R^d$ with $g = 0$ and $f = -1$. The exact (weak) solution is known and given by
\begin{align}\label{eq:fundamentalSolution}
\Phi_p (x) = \frac{p-1}{p d^{1/(p-1)}} (1 - |x|^{p/(p-1)}) .  
\end{align}
This solution for $d = 1$ is plotted for a range of $p$ in Fig.~\ref{fig:fundamentalSolution}. As $p$ increases from $p=2$, the cusp at the origin becomes sharper, highlighting that $\Phi_p (x)$ is necessarily a weak solution of \eqref{eq:pLaplacianDirichlet} for $p > 2$. (For $p \leq 2$, $\Phi_p \in C^{\infty} (B_1)$.) Note that we may write \eqref{eq:fundamentalSolution} in terms of $\delta$, where $p = 2 + \delta$ (i.e., the ordinary series), as
\begin{align*}
\Phi_{2 + \delta} (x) = \frac{1+\delta}{(2 + \delta) d^{\frac{1}{1+\delta}}} (1 - |x|^{\frac{2+\delta}{1 + \delta}}) .
\end{align*}

The map $\delta \mapsto \Phi_{2 + \delta} (x)$ is real analytic on $(-1 , \infty)$ for any $x \in B_1 (0)$. The radius of convergence about $\delta = 0$ with $d = 1$ is $|\delta| < 1$ ($1 < p < 3$). Letting $u_n = \frac{1}{n!} \partial_{\delta}^n \Phi_{2 + \delta} |_{\delta = 0}$, the first two terms in the Taylor series about $\delta = 0$ are $u_0 (x) = (1 - |x|^2) / 2d$ and 
\begin{align}\label{eq:pFSseries}
u_1 (x) = \dfrac{1}{4d} \left( (1 + \ln{d^2}) (1 - |x|^2) + |x|^2 \ln{|x|^2} \right) .
\end{align} 
Clearly, $u_0 (x) = (1 - |x|^2) / 2d$ satisfies $\Delta u_0 = -1$ and $u_0 |_{\partial B_1} = 0$. A straightforward calculation also demonstrates that $\Delta u_1 = -\diver{\ln\left|\nabla u_0\right|\nabla u_0}$ with $u_1 |_{\partial B_1} = 0$ (which has to be interpreted in a weak sense because $\R^d \ni x \mapsto |x|^2 \ln{|x|^2} \in [0,\infty)$ is not $C^2$ at $x = 0$). 

Reparameterizing \eqref{eq:fundamentalSolution} in terms of $p'$ reads
\begin{align}\label{eq:fundamentalSolution_prime}
\Phi_{p (p')} (x) = \frac{1}{p' d^{p'-1}} ( 1 - |x|^{p'} ) . 
\end{align}
Note the simpler expression \eqref{eq:fundamentalSolution_prime} takes compared to \eqref{eq:fundamentalSolution}. For this dual parameterization, the first couple of expansion coefficient functions $u_n' = \frac{1}{n!} \partial_{\delta}^n \Phi_{p(2 + \delta)} |_{\delta = 0}$ read $u_0' (x) = (1 - |x|^2) / 2d$ and
\begin{align}\label{eq:pprimeFSseries}
u_1' (x) = - \dfrac{1}{4d} \left( (1 + \ln{d^2}) (1 - |x|^2) + |x|^2 \ln{|x|^2} \right) .
\end{align}
Again, we may readily verify that $u_0'$ and $u_1'$ satisfy their respective PDEs from the hierarchy.

It is possible to compute the Taylor series coefficient functions $u_n$ and $u_n'$ corresponding to \eqref{eq:fundamentalSolution} and \eqref{eq:fundamentalSolution_prime}, respectively, exactly by using well-known relationships between power series coefficients of function compositions, products, etc. and their constituents (see, e.g., Fa\`{a} di Bruno's formula \eqref{eq:diBrunos}). The formula for $u_n$ is a bit cumbersome to write out explicitly for \eqref{eq:fundamentalSolution}. However, for \eqref{eq:fundamentalSolution_prime}, we have for $d = 1$ and $n \geq 1$,
\begin{align}\label{eq:pprimeFSseries_nth}
\scalemath{0.9}{
u_n' (x) = \frac{(-1)^n (1-|x|^2)}{2^{1+n}} - \frac{|x|^2}{2} \sum_{k = 0}^n \frac{ (-1)^{n - k} \ln^k (|x|) }{ 2^{n - k} k! } .
}
\end{align}
The expansion $\sum_{n \geq 0} (p' - 2)^n u_n' (x)$ converges for any $p' \in (0 , 4)$ for any $d \geq 1$. This, in particular, implies that $\Phi_p (x) = \sum_{n \geq 0} \left( p' (p) - 2 \right)^n u_n' (x)$ is valid for any $4/3 < p < \infty$ and any $x \in B_1 (0)$. 

\begin{figure}
  \begin{center}
\includegraphics[width=1\columnwidth]{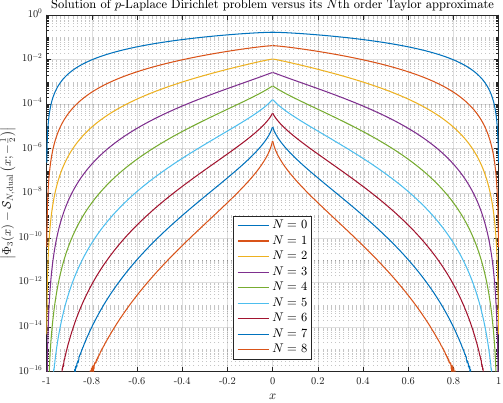}  
\caption{\label{fig:FSdiff}  \emph{Convergence of series expansion for an exact solution of the $p$-Laplacian Dirichlet problem.} The difference between the $N^{\mathrm{th}}$ partial sum $\mathcal{S}_{N,\mathrm{dual}} (x;-1/2)$ and the exact solution is shown for $p=3$. Adding additional terms clearly improves the agreement, with convergence being slowest near the cusp at the origin.}
\end{center}
\end{figure}

The zeroth-order functions $u_0$ are the same between ordinary and dual series, as expected, because both series are expanded about the same point $p = p' = 2$. The fact that the first order functions \eqref{eq:pFSseries} and \eqref{eq:pprimeFSseries} differ by a sign can be understood by this observation and comparing Eqs.~\eqref{eq:pLaplacian_Hierarchy_generalHomotopy} for $h (\delta) = \delta$ and $h(\delta) = -\delta / (1+ \delta)$. Similar to what was done in Section \ref{sec:BurgersFS}, we may use the exact expansion coefficients of $p \mapsto \Phi_p$ ($p' \mapsto \Phi_{p(p')}$) about $p = 2$ ($p' = 2$), respectively, to arbitrarily high order and plot the difference between $\Phi_p$ ($\Phi_{p'}$) and the $N^{\mathrm{th}}$ partial sum of the $k\in\{\mathrm{ordinary,dual}\}$ series, defined similar to Eq.~\eqref{eq:BurgersPartialSum} from the nonlinear advection case as
\begin{align}
\mathcal{S}_{N,k} (x;\delta) = \sum_{n=0}^{N}\delta^n u_n (x)\, .
\end{align}
Using \eqref{eq:pprimeFSseries_nth}, the results of this procedure are shown for the dual parameterization at $p=3$ ($p' = 3/2$) in Fig.~\ref{fig:FSdiff}. As will be discussed in more detail momentarily, additional terms to the series rapidly reduces the error. However, we see that convergence is slowest near the cusp at the origin and the error becomes clustered around this point as $N$ increases.

\begin{figure}
  \begin{center}
\includegraphics[width=1.0\columnwidth]{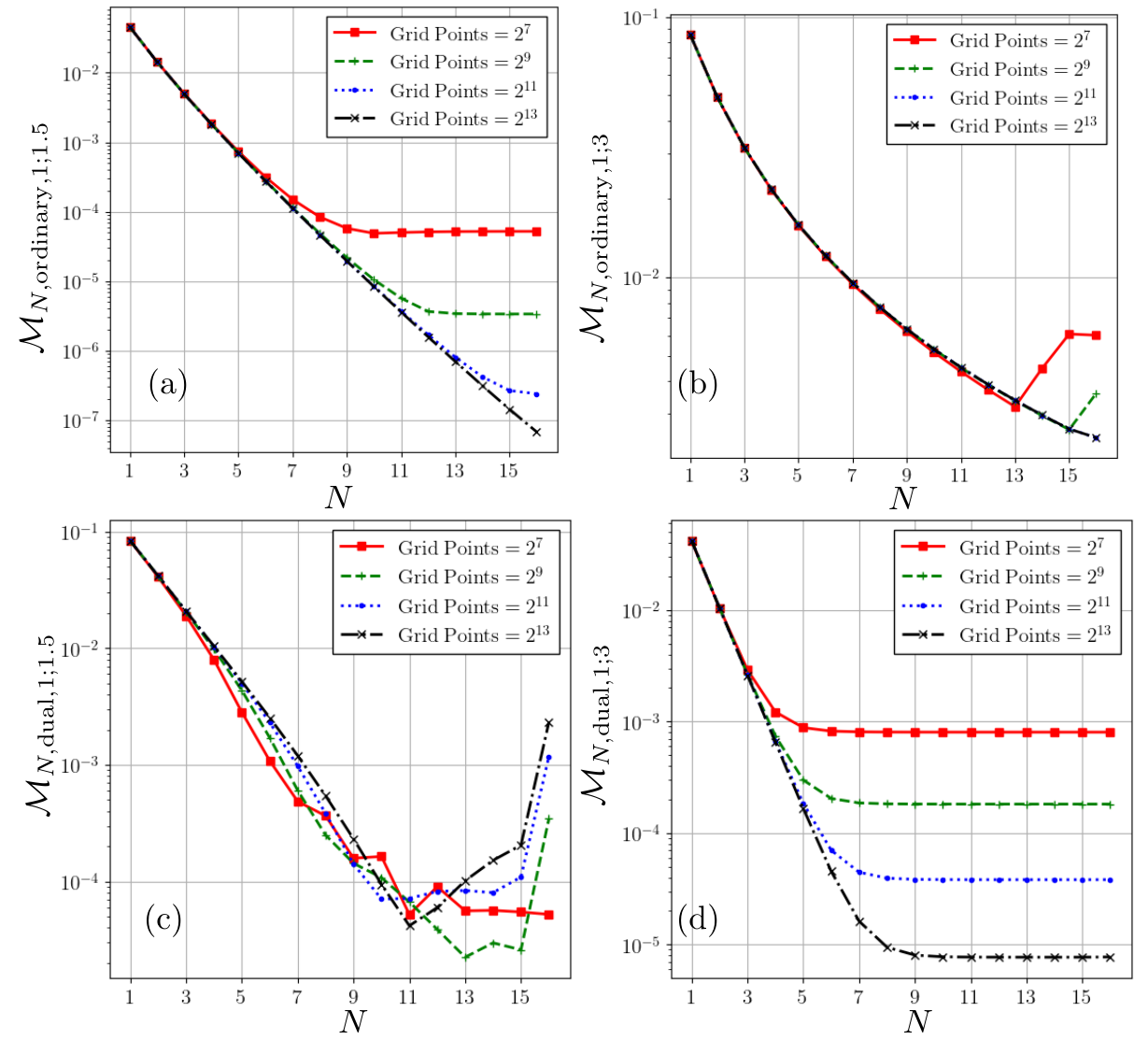}  
\caption{\label{fig:ErrorpLap}  \emph{Convergence with number of terms in the series $N$, $1$ dimension.} The $L_1$-norm error of the partial sums of the series expansions are shown vs.~the number of terms in the series for a range of $p$-values (left panels are $p=1.5$ and right panels are $p=3$), the ordinary and dual series (top panels are ordinary bottom panels are dual), and number of grid points used in the finite difference scheme (red $\square$ are $N_G=2^7$, green $+$ are $N_G=2^9$, blue $\cdot$ are $N_G=2^{11}$, and black $\times$ are $N_G=2^{13}$).  When the series is rapidly convergent ($p=1.5$ ordinary and $p=3$ dual), the finite differencing scheme leads to convergence ``plateaus" as $N$ increases.  }
\end{center}
\end{figure}

To numerically solve the hierarchy of linear Poisson-Dirichlet problems associated with \eqref{eq:HomotopyDirichlet_pLaplacian}, we employ a finite-difference scheme with standard second-order central differences for the Laplacian and first derivative operators \footnote{The standard finite-difference scheme operates on the solution pointwise, and so would appear inappropriate to treat PDEs involving the $p$-Laplacian which require a weak formulation. However, properties of our series methodology make the Dirichlet problem amenable to finite difference approaches provided the grid does not encompass points where a singularity occurs. The zeroth-order equation is a Poisson equation, which will have a smooth solution due to elliptic regularity (assuming a regular enough source term). All higher-order equations take the form $\nabla \cdot \nabla u_n=\nabla \cdot F_n$, which becomes the solution of a linear system $A \mathbf{a}=\mathbf{b}$ in the finite-difference formulation. For properly defined Laplacian and derivative operators of the same order and stencil, the effects of non-differentiability are regularized on both sides similar to how regularization occurs in a finite element method (FEM) formulation. A FEM method can also be applied to the Dirichlet problem; however, we have kept the finite difference formulation here due to its simplicity in scaling to the 2D problem.}. We define a scalar metric of convergence $\mathcal{M}_{N,k,q;p}\equiv \| \mathcal{S}_{N,k}\left(x; \delta_p\right)-\Phi_p (x) \|_{L^q (B_1)}$ between the truncated series solution and the known result Eq.~\eqref{eq:fundamentalSolution}.  Using the finite difference scheme, we will have an additional source of error due to the discretization with a finite number of grid points $N_G$ that will limit $\mathcal{M}_{N,k;p}$ to a nonzero value even as $N\to\infty$.  Typical behavior as we increase $N_G$ and $N$ is shown in Fig.~\ref{fig:ErrorpLap} for $p=1.5$ (left panels) and $p=3$ (right panels) for both the ordinary (top panels) and dual series (bottom panels) in $1$ dimension, $d=1$.  

In places where the series is rapidly convergent, e.g. panel (a) showing the ordinary series at $p=1.5$, we see that the effect of finite $N_G$ is to create a convergence ``plateau" in $\mathcal{M}_{N,k;p}$ as we increase $N$, with the convergence behavior before we reach the plateau being exponential in $N$.  When the series converges more slowly, as in panel (b) showing the ordinary series for $p=3$, the convergence is limited by the number of terms in the series rather than the finite differencing error at the same range of number of gridpoints $N_G$ (note the difference in scale on the $y$ axis).  The story is similar for the dual series, but the range of where the series is ``rapidly convergent," is different between the two.  Generally speaking, the ordinary series performs better for $p\in(1,2)$ and the dual series performs better for $p>2$.

\begin{figure}
  \begin{center}
\includegraphics[width=0.7\columnwidth]{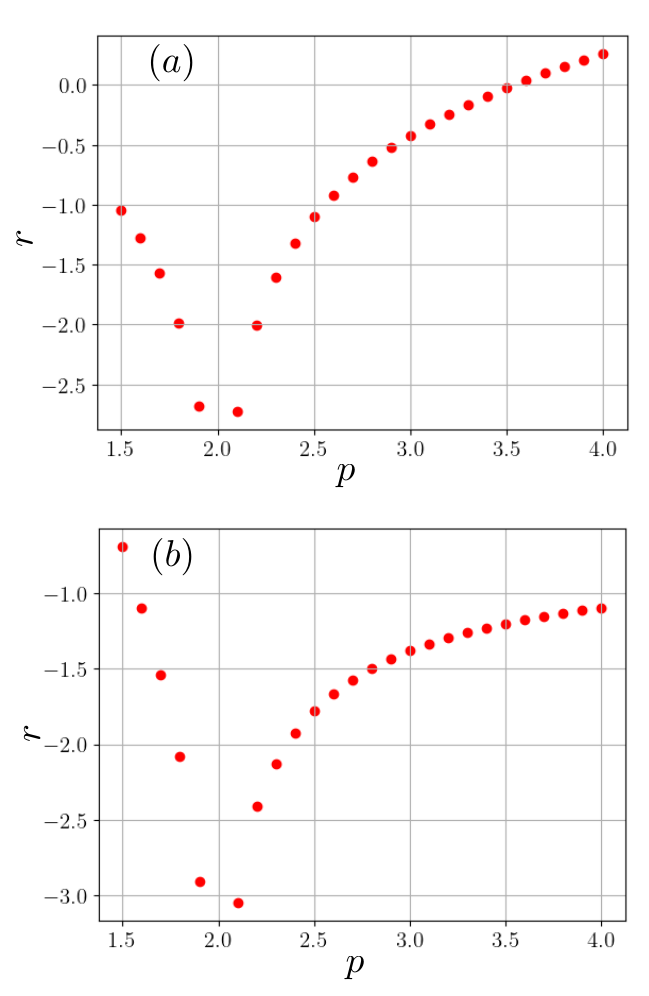}  
\caption{\label{fig:SlopepLap}  \emph{Rate of convergence vs.~$p$ for ordinary and dual series, $1$ dimension.} The convergence rate $r$ defined as in the main text for (a) regular series and (b) dual series. The points $p=2$ have been excluded, as the rate of convergence is zero but the error is also everywhere zero.}
\end{center}
\end{figure}

To characterize this rate of convergence we can estimate the slope of $\log |\mathcal{M}_{N,k,1;p}|$ vs.~$N$ to obtain an estimated rate $r$ as $|\mathcal{M}_{N,k,1;p}|\sim e^{rN}$, provided we are before the ``convergence plateau."  Practically, we can compute $r$ as $\sum_{n<N_{\mathrm{plateau}}} \log\left(\mathcal{M}_{n+1,k,1;p}/\mathcal{M}_{n,k,1;p}\right)/\sum_{n<n_{\mathrm{plateau}}} 1$, i.e. as the average slope of the convergence curve before the grid-driven plateau on a log-linear plot. These rates $r$ are plotted as a function of $p$ in Fig.~\ref{fig:SlopepLap} for the ordinary and dual series using $N_G=$ and $N_{\mathrm{plateau}}=5$.  We see that in the range $p\in(1.5,3)$ where $|\delta|\le 1$ for both the ordinary and the dual series, both approaches appear convergent in the sense that they have negative slope over this region.  We again stress that this metric only looks at the rate of convergence and not its magnitude; comparison with Fig.~\ref{fig:ErrorpLap} shows that the magnitude of error can be sizable at moderate $N$  even when the series appears to converge.  Convergence is excellent for both series near the point $p=2$, as expected.  As $p$ becomes large ($p\gtrsim 3$), the dual series continues to display convergence while the ordinary series fails to converge (slope becomes positive).

To investigate the role of dimensionality, we consider the two-dimensional version of Eq.~\eqref{eq:fundamentalSolution}, i.e. $d=2$.  By inspection, we can see that the solution separates into a radial function and a (trivial) angular function, and so would not be truly higher dimensional in those coordinates.  Hence, we consider a related non-separable problem on the unit square $\Omega=\left[-1,1\right]^2$:
\begin{align}
\label{eq:unitsquareKluge}\left\lbrace \begin{array}{ll}
\Delta_p (u) = -1 & \text{in } \Omega \\
u = \frac{p-1}{p d^{1/(p-1)}} (1 - |x|^{p/(p-1)})& \text{on }  \partial \Omega  ,
\end{array} \right.
\end{align}
Clearly, Eq.~\eqref{eq:fundamentalSolution} also solves this problem.  We again solve Eq.~\eqref{eq:unitsquareKluge} numerically using finite differences, applying the ordinary and dual series expansions. Given that our boundary data also depends on $p$, we have two choices of how to treat BCs: we can either set the boundary data of $u_0$ equal to the boundary data of the full problem or we can Taylor expand the boundary data and match BCs order-by-order. In what follows we follow the first choice.  We note that for this choice the numerically computed $u_n$ will not coincide with the Taylor expansion coefficient functions of $p \mapsto \Phi_p$ about $p = 2$.

\begin{figure}[t]
  \begin{center}
\includegraphics[width=1.0\columnwidth]{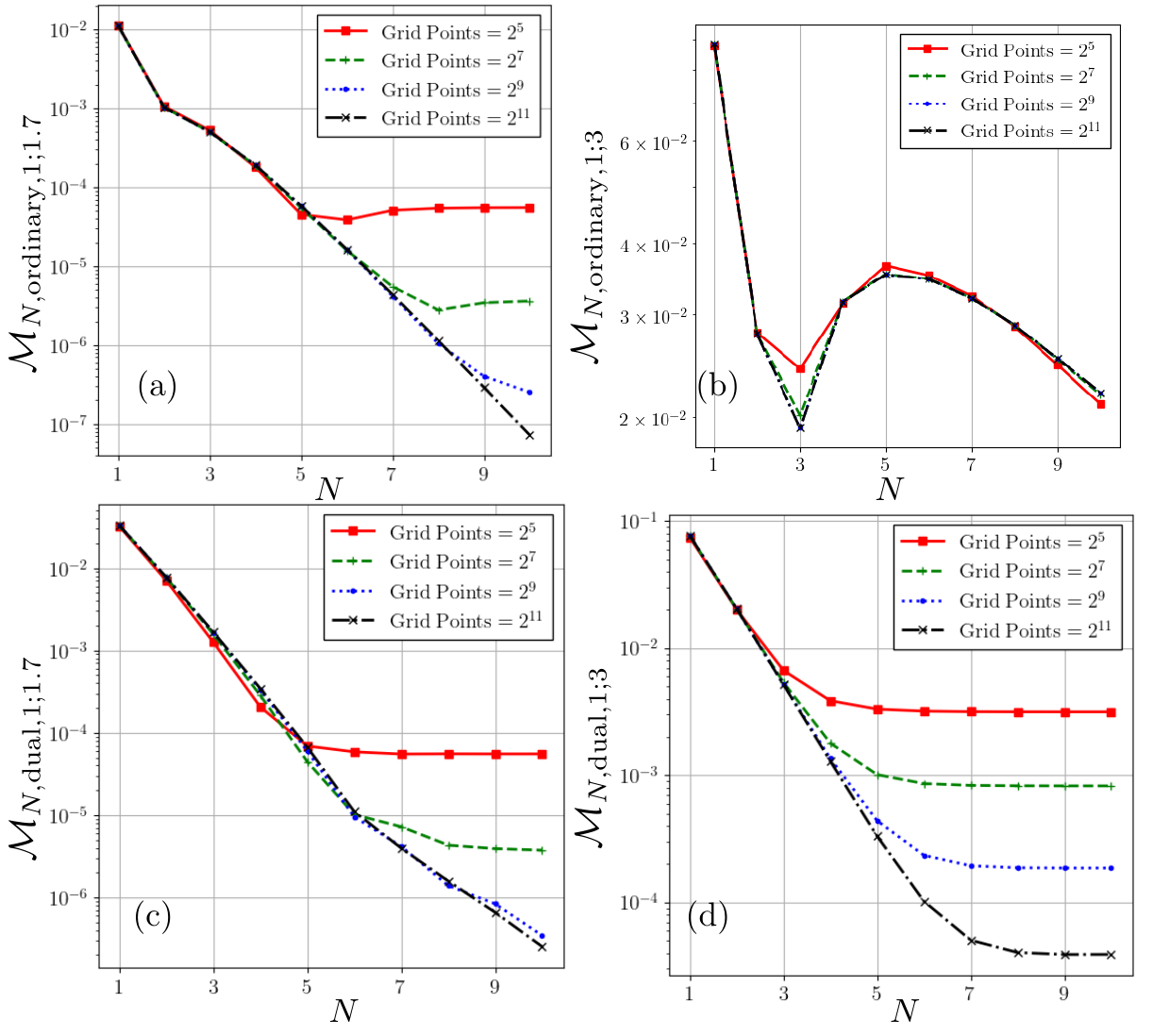}  
\caption{\label{fig:Error2D}  \emph{Convergence with number of terms in the series $N$, $2$ dimensions.} The $L_1$-norm error of the partial sums of the series expansions are shown vs.~the number of terms in the series for a range of $p$-values (left panels are $p=1.7$ and right panels are $p=3$), the ordinary and dual series (top panels are ordinary bottom panels are dual), and number of grid points used in the finite difference scheme.  Similar to the 1D case, the finite differencing scheme leads to convergence ``plateaus" as $N$ increases in regions where the series is rapidly convergent.}
\end{center}
\end{figure}

\begin{figure}[t]
  \begin{center}
\includegraphics[width=0.7\columnwidth]{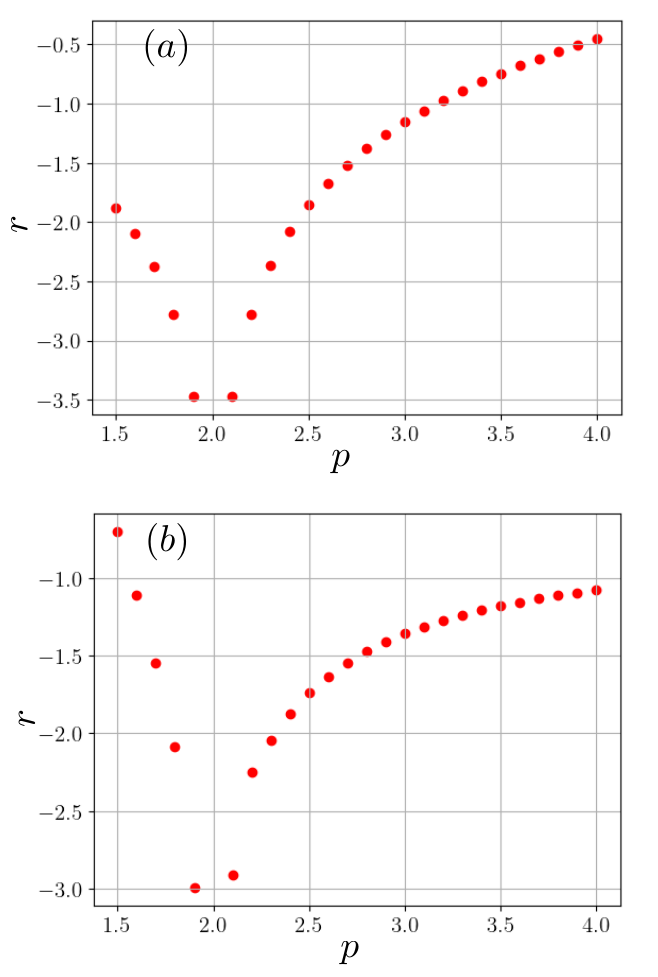}  
\caption{\label{fig:Slope2D}  \emph{Rate of convergence vs.~$p$ for ordinary and dual series, $2$ dimensions}. The convergence rate $r$ defined as in the main text for (a) regular series and (b) dual series for the two-dimensional Dirichlet problem Eq.~\eqref{eq:unitsquareKluge}.  The points $p=2$ have been excluded, as the rate of convergence is zero but the error is also everywhere zero.}
\end{center}
\end{figure}

A numerical comparison of the partial sums of the series approach with the known solution for the two-dimensional Dirichlet problem Eq.~\eqref{eq:unitsquareKluge}, analogous to Fig.~\ref{fig:ErrorpLap} for the one-dimensional case, are shown in Fig.~\ref{fig:Error2D}.  The behavior is very similar to the 1D case, with numerical ``plateaus" appearing due to finite grid size effects in regions where the series are rapidly convergent.  Unsurprisingly, numerical values of the error at fixed grid size are different between the 1D and 2D cases, with the ordinary series being less performant at $p=3$ in 2D than it was in 1D for the same grid spacing.  Similarly, by averaging the slope of the logarithm of the error, we can obtain an estimate of the convergence rates in the two-dimensional cases, analogous to the analysis presented in Fig.~\ref{fig:SlopepLap}.  The results for the two-dimensional case are shown in Fig.~\ref{fig:Slope2D}.  Qualitative results are again similar; the ordinary series performs best for $p<2$ while the dual series performs better for $p>2$.

\subsection{Convergence analysis - evolution equation}\label{sec:pLaplaceEvolution}

We now turn to the $p$-Laplacian evolution equation Eq.~\eqref{eq:pLaplacianEvolutionEqn}. Here, the case of Dirac-delta IC (i.e. the fundamental solution, also known in this context as the \emph{Barenblatt solution}~\cite{barenblatt1952self}) is known in closed-form and reads 
\begin{align}\label{eq:pLaplacianFundamentalSolution}
H_p (t,x) = t^{-k_p} \left( c_p - q_p |t^{-k_p/d} x|^{\frac{p}{p-1}} \right)_+^{\frac{p-1}{p-2}} ,
\end{align}
where
\begin{align}
\label{eq:kp} k_p = \left( p-2 + \frac{p}{d} \right)^{-1} \, ,\\
\label{eq:qp} q_p = \frac{p-2}{p} \left( \frac{k_p}{d} \right)^{\frac{1}{p-1}} \, ,
\end{align}
and the constant $c_p$ is determined by the condition
\begin{align}\label{eq:IntegralConstraint_FundamentalSoln}
\int_{\R^d} H_p (t,x) \dd x = 1 . 
\end{align}
In \eqref{eq:pLaplacianFundamentalSolution}, we make the assumption that $2d/(1+d) < p < \infty$. The lower bound here comes from the singular behavior of $k_p$ at $p = 2d/(1+d)$. Note that $1 \leq 2d / (1+d) < 2$ for any $d \in \N$. The subscript "$+$" in \eqref{eq:pLaplacianFundamentalSolution} indicates we are taking the positive part of the function in parentheses, i.e. $f_+ (x) := \max{\{0,f(x)\}}$. The constant $c_p$ does have a closed-form in terms of the Euler Beta function. The explicit form depends on whether $p > 2$ and $p < 2$. For $p > 2$,
\begin{align}
c_p^+ = q_p^{(p-2)k_p} \left( \frac{ p'(p) }{ |\mathbb{S}^{d-1}| \mathrm{B}\!\left( \frac{d}{p'(p)} , 1 + \frac{1}{2-p'(p)} \right) } \right)^{\lambda_p} ,
\end{align}
and, for $p < 2$,
\begin{align}
\scalemath{0.9}{
c_p^- = |q_p|^{(p-2)k_p} \left( \frac{ p'(p) }{ |\mathbb{S}^{d-1}| \mathrm{B}\!\left( \frac{d}{p'(p)} , \frac{1}{p'(p)-2} - \frac{d}{p'(p)} \right) } \right)^{\lambda_p} ,
}
\end{align}
where $|\mathbb{S}^{d-1}| = 2 \pi^{d/2} / \Gamma (d/2)$ is the surface area of the unit sphere in $\R^d$ and $\lambda_p = p(p-2)k_p / d(p-1)$.

The fundamental solution \eqref{eq:pLaplacianFundamentalSolution} has a number of intriguing properties~\cite{LEE2006389}. For one, it is a ``similarity solution" scaling as a universal function of $x$ and $t$ multiplied by another function of $t$, which comes from the scaling properties of the original equation.  The presence of the norm $|\cdot|$ highlights the general non-$C^2$ nature of the weak solution when $p > 2$.  Additionally, for $p>2$, $q_p>0$ and Eq.~\eqref{eq:pLaplacianFundamentalSolution} has compact support for any $t$. Contrariwise, when $p<2$, $q_p<0$ and the fundamental solution is non-compactly supported. A straightforward, if tedious, application of L'H$\hat{\mathrm{o}}$pital's rule on Eq.~\eqref{eq:pLaplacianFundamentalSolution} with $p \rightarrow 2^+$ and $p \rightarrow 2^-$ both yield $u_0 (t,x) = H_2(t,x)$, where $H_2$ is the classical heat kernel \eqref{eq:HeatKernel}, which is also non-compactly supported. A similar calculation, restricted to $1$ dimension for simplicity, yields that
\begin{widetext}
\begin{align}
u_1 (t,x) = \lim_{p \rightarrow 2^+} \frac{\dd}{\dd p} H_p (t,x) = \left( \ln{\sqrt{t}} - \frac{x^4 + 4 x^2 t \left( \ln{\left( \frac{16 \pi t^3}{x^2} \right)} - 1 \right) - 4 t^2 \left( \ln{(256 \pi^2)} + 2 \gamma - 3 \right)}{32 t^2}  \right) \frac{e^{-\frac{x^2}{4t}}}{\sqrt{4 \pi t}} ,
\end{align}
\end{widetext}
where $\gamma$ denotes Euler's constant, with numerical value $\gamma \simeq 0.577216$. One may likewise verify that the same limit is obtained from below, i.e. $u_1 (t,x) = \lim_{p \rightarrow 0^-} \frac{\dd}{\dd p} H_p (t,x)$. One can verify that satisfies $\partial_t u_1 - \partial_x^2 u_1 = \partial_x ( \ln{|\partial_x u_0|} \partial_x u_0 )$ with $u_0 (t,x) = H_2(t,x)$, the heat kernel \eqref{eq:HeatKernel}, in the appropriate weak sense. Note how $u_1$ is non-compactly supported, reflecting the fact that every equation in \eqref{eq:pLaplacian_Hierarchy_generalHomotopy} is a inhomogeneous heat equation and, hence, possesses infinite speed of propagation.

Solutions to the $p$-Laplacian evolution equation for $p > 2$ with compactly supported initial data remain compactly supported, a property known as finite speed of propagation.  This can be starkly contrasted with the standard heat equation (i.e., $p=2$), where the action of evolution is to convolve the initial distribution with a Gaussian heat kernel and so solutions are non-compactly supported at any time $t>0$, a property shared with solutions of the $p$-Laplacian evolution equation with $p\le 2$.  Since the compact data becomes non-compact instantaneously, this phenomenon can be referred to as infinite speed of propagation.  Interestingly, the hierarchy \eqref{eq:pLaplacian_Hierarchy_generalHomotopy} must capture \emph{both} phenomena within the same set of functions $\{u_n(x)\}$, where it is known that the functions are non-compactly supported because they are solutions to inhomogeneous heat equations. This implies that sequential contributions to the series outside the compact region of support must exhibit some cancellation behavior. We will see this behavior arise numerically in the analysis that follows.


To numerically solve the ordinary and dual series hierarchies for the $p$-Laplacian evolution equation in $1$ dimension we employ a finite element method (FEM) based on piecewise linear ``hat" functions.  Each hat function $\phi_i(x)$ takes the value $1$ at the point $x_i=-L+i\Delta x$, decreases linearly to zero at neighboring points $x_{i-1}$ and $x_{i+1}$, and is zero elsewhere in the domain.  $L$ is chosen large enough to capture the support of the needed functions to the desired numerical accuracy and BCs set accordingly.  We then represent a function in this basis as $u(t,x) = \sum_i a_i (t) \phi_i(x)$ and convert the weak formulation of our equations to linear systems of ODEs for the $a_i(t)$ in a Galerkin scheme.  Defining the matrices $A_{ij}=\langle \phi_i(x),\phi_j(x)\rangle$, $B_{ij}=\langle \nabla \phi_i(x),\nabla \phi_j(x)\rangle$, and $C_{ij}=\langle \nabla \phi_i(x), \phi_j(x)\rangle$, with $\langle \cdot,\cdot\rangle$ being the $L^2$-inner product, we convert the linear, driven diffusion equation
\begin{align}
\partial_t u(t,x)&=\Delta u(t,x) + \diver{F(t,x)}\, ,
\end{align}
into the ODEs
\begin{align}\label{eq:GalerkinFormulation}
A \dot{\mathbf{a}} &= B \mathbf{a}(t) - C \mathbf{f}(t)
\end{align}
in which $\mathbf{a}$ is the vector with $a_i(t)$ as components and the vector $\mathbf{f}(t)$ is defined through the Galerkin expansion of the forcing,
\begin{align}
F(t,x)&=\sum_{i} f_i (t) \phi_i(x)\, .
\end{align}
We solve the formulation Eq.~\eqref{eq:GalerkinFormulation} using an implicit Euler method 
\begin{align}
\left[A-\Delta t B\right]\mathbf{a}(t+\Delta t)&=A\mathbf{a}(t)-\Delta t C\mathbf{f}(t+\Delta t)\, ,
\end{align}
with small timesteps $\Delta t$. Because of our choice of basis functions, we have that $u(x_i,t)=a_i(t)$ for a well-converged solution, and we use this representation with simple numerical quadrature when computing metrics. Similarly, while the forcing functions appearing in our equations are functionals of functions obtained as outputs of an FEM process (see Eq.~\eqref{eq:un_diff}) and so should strictly be obtained by $L^2$ projection of the functional onto the Galerkin basis, we have found that evaluating those functionals pointwise does not introduce appreciable additional error for our parameters and is significantly faster.  We do initialize with $[a_0(t=0)]_i = \left[A^{-1}\right]_{i,i_0}$, where $i_0$ is the grid point at $x=0$, which is the appropriate $L^2$ projection of the delta function IC onto the Galerkin basis, and perform a similar projection when computing gradients of functions obtained from the FEM prcoess.

\begin{figure*}
  \begin{center}
\includegraphics[width=2\columnwidth]{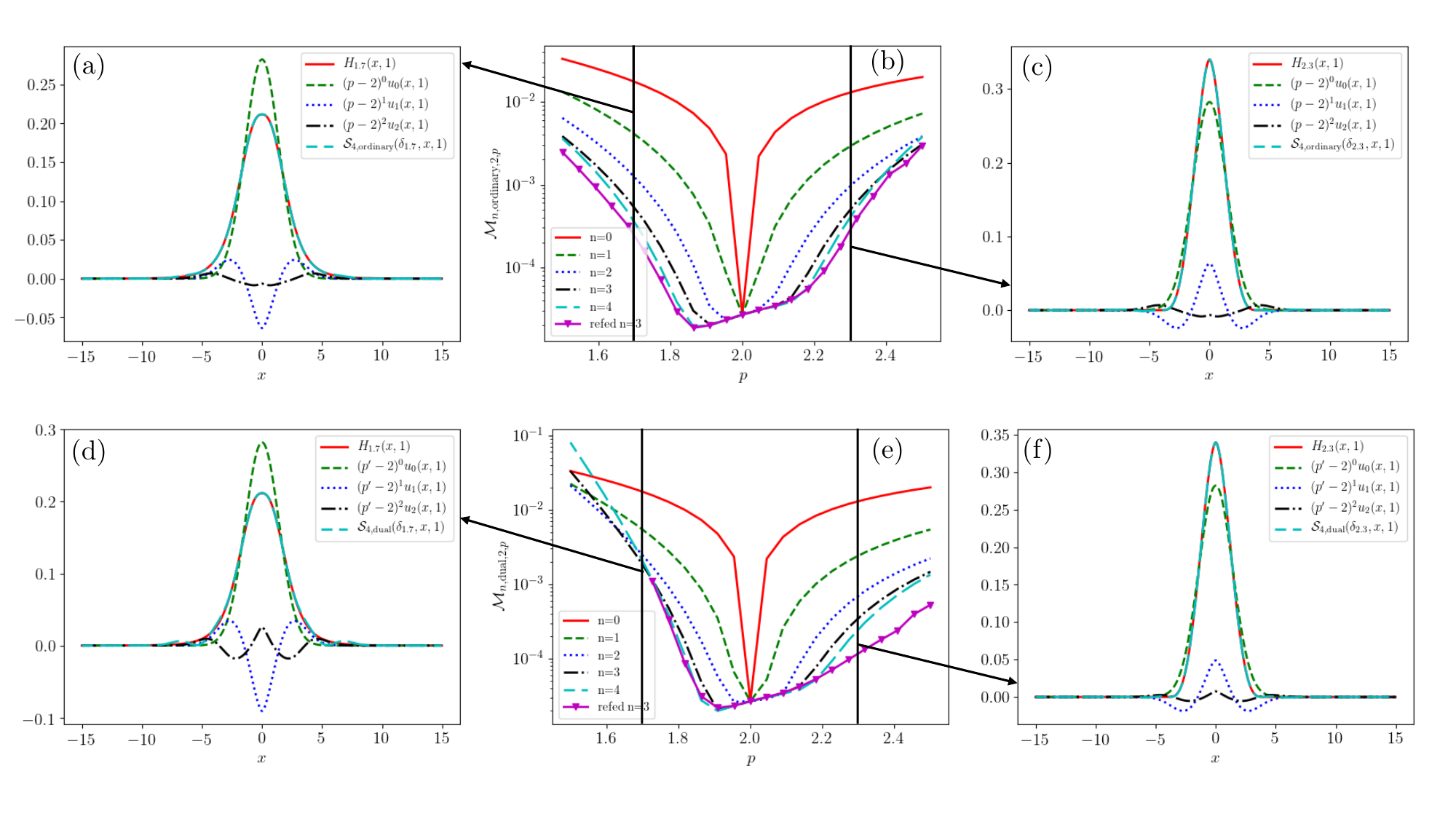}  
\caption{\label{fig:Barenblatt}  \emph{Behavior of series solutions for the $p$-Laplacian evolution equation with delta-function IC.}. Numerical solutions of our series approach are compared against the exact solution $H_p(t,x)$ from Eq.~\eqref{eq:pLaplacianFundamentalSolution} at $t=1$ for a range of $p$ using the ordinary homotopy (top row) and the dual series (bottom row).  Panels (a) and (d) compare the exact solution with the series partial sums to fourth order along with some of the terms in the partial series in the regime $p<2$, with the top being the ordinary series and the bottom the dual series.  Panels (c) and (f) are analogous for the regime $p>2$.  The center panels (b) and (e) show the behavior of the $L_2$ norm residuals between the exact solution and the partial series up to order 4 for a range of $p$ together with the solution from the refeeding process up to order $n=3$.  The vertical solid lines indicate the $p$ where the solutions are shown in the left and right panels.  As was the case for the Dirichlet problem, the ordinary series generally performs better for $p<2$ and the dual series performs better for $p>2$.}
\end{center}
\end{figure*}

Fig.~\ref{fig:Barenblatt} shows the results of our FEM methodology with $\Delta x \sim 0.02$ and $\Delta t\sim 0.01$ compared with the closed-form solution $H_p(t,x)$ from Eq.~\eqref{eq:pLaplacianFundamentalSolution}.  The top row of panels shows the ordinary series, and the bottom row the dual series.  The center panels display the residuals of the partial sums through fourth order with the exact fundamental solution at time $t=1$.  Qualitatively, we see similar behavior to what was seen for the Dirichlet case: there is a region near $p=2$ where both series perform well and errors are monotonically decreasing.  We have also implemented the refeeding process described in Sec.~\ref{sec:refeeding} up to order $n=3$ and refeeding every timestep, and the errors are shown alongside the errors without refeeding.  The solutions resulting from refeeding perform at least as well as the full series expansions up to order $n=4$ in spite of only being carried out to order $n=3$, except for the case of the dual series with $p\lesssim 1.7$ where refeeding did not converge.  We can see that the series solution with the given simulation parameters is likely unreliable in this region, as the $n=4$ residual becomes larger than the $n=3$ residual.  We note that the interplay between convergence parameters (time step, spatial grid size, domain length, etc.) may be expected to be especially delicate for the chosen problem on account of the singular initial condition, and other choices of parameters or FEM methodology may be expected to improve the convergence.  We leave more detailed analysis of such convergence for future work.  Generally speaking, the ordinary series performs better further to the left of $p=2$, while the dual series performs better further to the right of $p=2$. For $p<2$, solutions are non-compactly supported and shorter and wider than the Gaussian solution for $p=2$.  Contrariwise, for $p>2$ the exact solution is compactly supported and so is narrower and taller than the Gaussian at $p=2$. The left and right panels display how the functions $u_n$ are weighted in the series to reproduce this behavior well in the partial sum.  

\section{Discussion and outlook}
\label{sec:DaO}

We considered nonlinear advection-diffusion PDEs which are augmented by an additional parameter $\delta$ such that $\delta=0$ produces linear advection-diffusion, and showed how a series expansion in $\delta$ matched order-by-order obtains the solution of the nonlinear PDE through solving a hierarchy of linear, forced PDEs and summing their solutions.  We studied the case of nonlinear advection through Burgers' equation in detail, and proved that a novel linear deformation encompassing linear advection-diffusion and Burgers' equation is analytic in $\delta$ and has infinite radius of convergence.  We further showed how the hierarchy of linear equations in this case can be solved efficiently using spectral techniques and long times can be accessed by iteratively restarting the series expansion in a time-stepping approach we call refeeding.  Our approach was shown to reproduce known exact results as well as behavior characteristic of Burgers' turbulence in a forced scenario with no known closed-form solution.  We then considered nonlinear diffusion generated by the $p$-Laplacian operator in both static and dynamic contexts.  Here, we analyzed two different deformations of the nonlinear equation connecting it to linear diffusion that we called the ordinary and dual, and analyzed performance against known exact solutions for statics in $1$ and $2$ dimensions and dynamics in $1+1$ dimensions.  The dual case involves a nonlinear function of the deformation parameter $\delta$, and enables efficient simulations outside of a naive perturbative regime.

Our work sets a rigorous foundation for using series methods to study nonlinear PDEs, and opens many further avenues for research.  A clear next step is to apply series methods to higher-dimensional PDEs that include both nonlinear advection and nonlinear diffusion, such as those arising from LES models, or to apply the method to directly treat the nonlinear advection in NS. While the focus of the present work was on problems with known exact solutions for benchmarking purposes, it will be interesting to see how series methods perform against DNS in challenging regimes of fluid flow. Given that our methods reduce the solution of a nonlinear PDE to a system of linear PDEs, future work will also explore means of leveraging the linear representations for quantum-enhanced computations.

While our focus has been on advection-diffusion PDEs due to their central role in pure and applied science, we note that other PDEs can be treated in a similar fashion, and we hope our work provides motivation for researchers in other fields to apply similar series approaches. Further, we note that our methods can also be applied to the solution of systems of nonlinear ODEs, for which analytic dependence of solutions on parameters of the equation is well established~\cite{Coddington1956theory}. This contrasts with the corresponding theory for PDEs, for which results are available in the linear case \cite{HormanderALPDOII}, but the nonlinear case appears much less developed.


\section{Acknowledgements}

We would like to thank Frank Graziani, Jeanie Qi, and Ilon Joseph for valuable discussions. This material is based upon work supported by the U.S. Department of Energy, Office of Science, Office of Fusion Energy Sciences, under Award Number(s) DE-SC0025203.

\appendix
\section{Series expansion of a generalized deformation of Burgers' equation}
\label{app:GeneralBE}

In this appendix we will derive the equation for $u_k$, $k \geq 1$, in the hierarchy of equations \eqref{eq:BurgerHomotopyHierarchy} resulting from a general homotopy connecting linear advection-diffusion and Burgers' equation. The equation for the $k^{\mathrm{th}}$-order $u_k$ in \eqref{eq:BurgerHomotopyHierarchy} reads
\begin{align*}
\partial_t u_k + \frac{1}{k!} \partial_x \frac{\dd^k}{\dd \delta^k} h (u (\delta) , \delta) \Big|_{\delta = 0}  - \nu \partial_x^2 u_k = 0 .
\end{align*}
To derive a general expression for the $k^{\mathrm{th}}$ derivative of $h(u(\delta),\delta)$ with respect to $\delta$, consider the Taylor expansion of $h$ about the point $(u_0 , 0)$:
\begin{align*}
h (u (\delta) , \delta) = \sum_{n_1 , n_2 \geq 0} \frac{(u(\delta) - u_0)^{n_1} \delta^{n_2} }{n_1! n_2!} h^{(n_1,n_2)}_0 ,
\end{align*}
where $h^{(n_1,n_2)}_0 := ( \partial_1^{n_1} \partial_2^{n_2} h ) (u_0 , 0)$. Using the binomial theorem applied to the $k^{\mathrm{th}}$ derivative of a product, we find
\begin{align*}
\frac{\dd^k}{\dd \delta^k} h (u (\delta) , \delta) & = \sum_{n_1 , n_2 \geq 0} \frac{h^{(n_1,n_2)}_0}{n_1! n_2!} \frac{\dd^k}{\dd \delta^k} \left[ (u(\delta) - u_0)^{n_1} \delta^{n_2} \right] \\
& = \sum_{n_1 , n_2 \geq 0} \frac{h^{(n_1,n_2)}_0}{n_1! n_2!} \sum_{m = 0}^k \binom{k}{m}  \\
& \hspace{5mm} \times \left( \frac{\dd^{k-m}}{\dd \delta^{k-m}}(u(\delta) - u_0)^{n_1} \right) \left( \frac{\dd^m}{\dd \delta^m} \delta^{n_2} \right) .
\end{align*}
Now,
\begin{align*}
\frac{\dd^m}{\dd \delta^m} \delta^{n_2} = \left\lbrace \begin{array}{ll}
\dfrac{n_2!}{(n_2 - m)!} \delta^{n_2 - m} & \text{if } m \leq n_2 , \\[5pt]
0 & \text{otherwise} .
\end{array} \right.
\end{align*}
Hence,
\begin{align*}
\frac{\dd^k}{\dd \delta^k} h (u (\delta) , \delta) \Big|_{\delta = 0} & = \sum_{n_1 \geq 0} \sum_{m = 0}^k \frac{h_0^{(n_1,m)}}{n_1!} \binom{k}{m} \\
& \hspace{5mm} \times \frac{\dd^{k-m}}{\dd \delta^{k-m}} (u(\delta) - u_0)^{n_1} \Big|_{\delta = 0} .
\end{align*}

The next step is to apply Fa\`{a} di Bruno's formula. Let $f,g : \R \rightarrow \R$ be any smooth functions. Then, for $n \geq 1$ integer, Fa\`{a} di Bruno's formula reads
\begin{align}\label{eq:diBrunos}
& \nonumber \frac{1}{n!} \frac{\dd^n}{\dd x^n} f(g(x)) \\
& \hspace{5mm} = \sum_{\mathbf{p} \in S_n} f^{(\| \mathbf{p} \|_1)} (g(x)) \prod_{j=1}^n \frac{1}{p_j!} \left( \frac{g^{(j)} (x)}{j!} \right)^{p_j} ,
\end{align}
where $\mathbf{p} = (p_1 , \cdots , p_n) \in \Z_+^n$, $\| \mathbf{p} \|_1 = \sum_j p_j$, and
\begin{align}\label{eq:S_n}
S_n = \{ \mathbf{p} \in \Z_+^n ~ : ~ \sum_{j=1}^n j p_j = n \} 
\end{align}
is the set of all integer partitions of $n$ (counting multiplicity). Applying \eqref{eq:diBrunos}, we find
\begin{align*}
&\frac{1}{(k-m)!} \frac{\dd^{k-m}}{\dd \delta^{k-m}} (u(\delta) - u_0)^{n_1} = \\
\nonumber &\sum_{\substack{\mathbf{p} \in S_{k-m} \\ \| \mathbf{p} \|_1 \leq n_1}} \frac{n_1! (u(\delta) - u_0)^{n_1 - \| \mathbf{p} \|_1}}{(n_1 - \| \mathbf{p} \|_1)!} \prod_{j=1}^{k-m} \frac{1}{p_j!} \left( \frac{u^{(j)} (\delta)}{j!} \right)^{p_j} .
\end{align*}
Therefore,
\begin{align*}
\frac{1}{n_1 ! (k-m)!} \frac{\dd^{k-m}}{\dd \delta^{k-m}} (u(\delta) - u_0)^{n_1} \Big|_{\delta = 0} = \sum_{\substack{\mathbf{p} \in S_{k-m} \\ \| \mathbf{p} \|_1 = n_1}} \prod_{j=1}^{k-m} \frac{u_j^{p_j}}{p_j!} .
\end{align*}
Hence, the desired formula for the $k^{\mathrm{th}}$ derivative reads
\begin{align}\label{eq:kthDerivativeHomotopy}
& \nonumber \frac{1}{k!} \frac{\dd^k}{\dd \delta^k} h (u (\delta) , \delta) \Big|_{\delta = 0} \\
& \hspace{5mm} = \sum_{m = 0}^k \sum_{\mathbf{p} \in S_{k-m}} \frac{h_0^{(\|\mathbf{p}\|_1,m)}}{m!} \prod_{j=1}^{k-m} \frac{u_j^{p_j}}{p_j!} .
\end{align}

Let's now identify which terms in \eqref{eq:kthDerivativeHomotopy} depend on the current order $u_k$. The current order $u_k$ appears in \eqref{eq:kthDerivativeHomotopy} only when $m = 0$. When $m = 0$, the first term on the right hand side of \eqref{eq:kthDerivativeHomotopy} simplifies to
\begin{align*}
\nonumber &\sum_{\mathbf{p} \in S_k} h_0^{(\mathbf{p}\|_1,0)} \prod_{j=1}^k \frac{u_j^{p_j}}{p_j!} \\
& = \sum_{\substack{ \mathbf{p} \in S_k \\ \| \mathbf{p} \|_1 = 1}} h_0^{(\mathbf{p}\|_1,0)} \prod_{j=1}^k \frac{u_j^{p_j}}{p_j!} + \sum_{\substack{ \mathbf{p} \in S_k \\ \| \mathbf{p} \|_1 > 1}} h_0^{(\mathbf{p}\|_1,0)} \prod_{j=1}^k \frac{u_j^{p_j}}{p_j!} .
\end{align*}
The first term on the right hand side of the previous simplifies to
\begin{align*}
\sum_{\substack{ \mathbf{p} \in S_k \\ \| \mathbf{p} \|_1 = 1}} h_0^{(\mathbf{p}\|_1,0)} \prod_{j=1}^k \frac{u_j^{p_j}}{p_j!} = h_0^{(1,0)} u_k
\end{align*}
because the only time $\| \mathbf{p} \|_1 = 1$ is when $p_k = 1$ and $p_j = 0$ for all $1 \leq j \leq k-1$, i.e. when $\mathbf{p} = (0 , \cdots , 0 , 1)$. Hence, there is only one term in \eqref{eq:kthDerivativeHomotopy} which depends on $u_k$, the current order, and that term reads $h_0^{(1,0)} u_k$ (this is consistent with the pattern emerging in \eqref{eq:BurgerHomotopyHierarchy}). Introduce the notation
\begin{align}\label{eq:InhomogeneityGeneralHomotopy}
\nonumber F_k (u_0 , \cdots , u_{k-1}) & := \sum_{\substack{ \mathbf{p} \in S_k \\ \| \mathbf{p} \|_1 > 1}} h_0^{(\mathbf{p}\|_1,0)} \prod_{j=1}^k \frac{u_j^{p_j}}{p_j!}   \\
& + \sum_{m = 1}^k \sum_{\mathbf{p} \in S_{k-m}} \frac{h_0^{(\mathbf{p}\|_1,m)}}{m!} \prod_{j=1}^{k-m} \frac{u_j^{p_j}}{p_j!} ,
\end{align}
with the understanding that when $m = k$ the term contributing to the summation is $h_0^{(0,k)} / k!$. With this notation, the equation for $u_k$ in \eqref{eq:BurgerHomotopyHierarchy} reads
\begin{align}\label{eq:Equation_for_uk_GeneralHomotopy}
\partial_t u_k + \partial_x \left( h_0^{(1,0)} u_k \right)  - \nu \partial_x^2 u_k = - \partial_x F_k (u_0 , \cdots , u_{k-1}) . 
\end{align}

Let's now specialize to the linear homotopy case, i.e. $H(u,\delta) = (1-\delta) v u + \delta u^2 / 2$, and look at the form of \eqref{eq:InhomogeneityGeneralHomotopy}. This case will also serve as a sanity check on formula \eqref{eq:InhomogeneityGeneralHomotopy}. Let's first collect the non-zero partials of $H$:
\begin{align*}
\left\lbrace \begin{array}{l}
h_0^{(1,0)} = v \\[5pt]
h_0^{(0,1)} = \dfrac{1}{2} u_0^2 - v u_0 \\[5pt]
h_0^{(1,1)} = u_0 - v \\[5pt]
h_0^{(2,1)} = 1 .
\end{array} \right.
\end{align*}
From these relations, we may immediately note that the second term on the right hand side of \eqref{eq:InhomogeneityGeneralHomotopy} is zero. Moreover, $F_1 (u_0) = - v u_0 + u_0^2 / 2 ,$ and, for $k \geq 2$,
\begin{align*}
&F_k (u_0 , \cdots , u_{k-1}) \\
& = \sum_{m = 1}^k \sum_{\mathbf{p} \in S_{k-m}} \frac{1}{m!} h_0^{(\mathbf{p}\|_1,m)} \prod_{j=1}^{k-m} \frac{u_j^{p_j}}{p_j!} \\
& = \sum_{\mathbf{p} \in S_{k-1}} h_0^{(\| \mathbf{p} \|_1,1)} \prod_{j=1}^{k-1} \frac{u_j^{p_j}}{p_j!} \\
& = (u_0 - v) \sum_{\substack{ \mathbf{p} \in S_{k-1} \\ \|\mathbf{p}\|_1 = 1}} \prod_{j=1}^{k-1} \frac{u_j^{p_j}}{p_j!} + \sum_{\substack{ \mathbf{p} \in S_{k-1} \\ \|\mathbf{p}\|_1 = 2}} \prod_{j=1}^{k-1} \frac{u_j^{p_j}}{p_j!} \\
& = (u_0 - v) u_{k-1} + \sum_{\substack{ \mathbf{p} \in S_{k-1} \\ \|\mathbf{p}\|_1 = 2}} \prod_{j=1}^{k-1} \frac{u_j^{p_j}}{p_j!} . 
\end{align*}
To continue simplifying, we need to understand the set $\{ \mathbf{p} \in S_n ~ : ~ \| \mathbf{p} \|_1 = 2 \}$. This set may be described in words as the number of ways to partition the integer $n$ into two smaller integers (counting multiplicity). Each vector $\mathbf{p} \in S_n$ with $\| \mathbf{p} \|_1 = 2$ may be expressed as $\mathbf{p} = \mathbf{e}_j + \mathbf{e}_{n-j}$ for $j = 1 , 2 , \cdots , \lfloor n/2 \rfloor$, where $\mathbf{e}_j \in \Z^n_+$ is the vector with $1$ in the $j$th component and zeros in all other components. Hence, for $k \geq 2$,
\begin{align*}
&F_k (u_0 , \cdots , u_{k-1}) = (u_0 - v) u_{k-1}  \\
&+\sum_{m = 1}^{\lfloor (k-1)/2 \rfloor} \prod_{j=1}^{k-1} \frac{u_j^{ \delta_{jm} + \delta_{j(k-1-m)} } }{(\delta_{jm} + \delta_{j(k-1-m)})!} . 
\end{align*}
This further simplifies to
\begin{align}
\label{eq:InhomogeneityLinearHomotopy} &F_k (u_0 , \cdots , u_{k-1}) = \\
\nonumber &\left\lbrace \begin{array}{ll}
(u_0 - v) u_{k-1} + \dfrac{1}{2} u_{(k-1)/2}^2 + \sum\limits_{m = 1}^{(k-3)/2} u_m u_{k-1-m}  & k \text{ odd,} \\[10pt]
(u_0 - v) u_{k-1} + \sum\limits_{m = 1}^{(k-2)/2} u_m u_{k-1-m} & k \text{ even.}
\end{array} \right.
\end{align}
The first few forms of $F_k$ in \eqref{eq:InhomogeneityLinearHomotopy} read
\begin{align*}
F_1 (u_0) & = - v u_0 + \frac{1}{2} u_0^2 , \\
F_2 (u_0,u_1) & = (u_0 - v) u_1 , \\
F_3 (u_0,u_1,u_2) & = (u_0 - v) u_2 + \frac{1}{2} u_1^2 , \\
F_4 (u_0,u_1,u_2,u_3) & = (u_0 - v)u_3 + u_1 u_2 , \\
F_5 (u_0,u_1,u_2,u_3,u_4) & = (u_0 - v) u_4 + \frac{1}{2} u_2^2 + u_1 u_3 . 
\end{align*}
Hence, the first $6$ orders in the hierarchy for the linear homotopy case (i.e., Eq.~\eqref{eq:BurgerLinearHomotopy_dimensionless}) read
\begin{align}\label{eq:LinearHomotopyHierarchy}
\scalemath{0.95}{ \left\lbrace \begin{array}{l}
\partial_t u_0 + \mathcal{L}_0 u_0 = 0 \\[5pt]
\partial_t u_1 + \mathcal{L}_0 u_1 = (v - u_0) \partial_x u_0 \\[5pt]
\partial_t u_2 + \mathcal{L}_0 u_2 = - \partial_x \left( (u_0 - v) u_1 \right) \\[5pt]
\partial_t u_3 + \mathcal{L}_0 u_3 = - \partial_x \left( (u_0 - v) u_2 + \dfrac{1}{2} u_1^2 \right) \\[5pt]
\partial_t u_4 + \mathcal{L}_0 u_4 = - \partial_x \left( (u_0 - v)u_3 + u_1 u_2 \right) \\[5pt]
\partial_t u_5 + \mathcal{L}_0 u_5 = - \partial_x \left( (u_0 - v) u_4 + \dfrac{1}{2} u_2^2 + u_1 u_3 \right) \\
\hspace{2cm} \vdots 
\end{array}  \right. 
}
\end{align} 
Formulas \eqref{eq:InhomogeneityLinearHomotopy}-\eqref{eq:LinearHomotopyHierarchy} can be additionally verified by directly differentiating the linear homotopy function $H(u (\delta),\delta) = (1-\delta) v u (\delta) + \delta u (\delta)^2 / 2$ with respect to $\delta$ and evaluating at $\delta = 0$.

\section{Series expansion of a generalized deformation of the $p$-Laplacian}
\label{app:seriespLap}
Here, we detail a method to obtain the PDE that all orders $u_n$ in the series expansion of the homotopy $p$-Laplacian evolution equation Eq.~\eqref{eq:HomotopyEvolution_pLaplacian} must satisfy for a general homotopy $h (\delta)$, analogous to Appendix \ref{app:GeneralBE} for the homotopy Burgers' equation. First, we set the notation that $f(\delta) = |\nabla u(\delta)|^{h (\delta)} \nabla u (\delta)$ and $u(\delta) \equiv u(t,x;\delta)$ is the solution to \eqref{eq:HomotopyEvolution_pLaplacian}.  We now investigate the $f^{(n)}(\delta)$ with an eye towards taking $\delta\to 0$.  We find
\begin{align}
\frac{\dd^n}{\dd\delta^n} f(\delta)&=\sum_{k=0}^{n}\binom{n}{k}\left(\frac{\dd^k}{\dd\delta^k}\left|\nabla u(\delta)\right|^{h (\delta)}\right)\left(\frac{\dd^{n-k}}{\dd\delta^{n-k}} \nabla u(\delta)\right)\, .
\end{align}
Using the expansion $u(\delta)=\sum_{n\ge 0} \delta^n u_n(\delta)$, we find
\begin{align}
\frac{\dd^{n-k}}{\dd\delta^{n-k}} \nabla u(\delta)\Big|_{\delta=0}&=(n-k)! \nabla u_{n-k}\, ,
\end{align}
and so
\begin{align}
\frac{\dd^n}{\dd\delta^n}f(\delta)&=\sum_{k=0}^{n}\frac{n!}{k!}\left(\frac{\dd^k}{\dd\delta^k}\left|\nabla u(\delta)\right|^{h (\delta)}\right)\nabla u_{n-k}\, , \\
\Rightarrow f(\delta)&=\sum_{n\ge 0} \delta^n \sum_{k=0}^{n}\frac{1}{k!}\left(\frac{\dd^k}{\dd\delta^k}\left|\nabla u(\delta)\right|^{h (\delta)}\right)\nabla u_{n-k}\, .
\end{align}
Now, using
\begin{align}
\left|\nabla u(\delta)\right|^{h (\delta)}&=\exp\left(h (\delta) \ln \left|\nabla u(\delta)\right|\right)\, ,\\
&=\sum_{m=0}^{\infty} \frac{1}{m!} h^m(\delta)\left(\ln \left|\nabla u(\delta)\right|\right)^m\, ,
\end{align}
we find
\begin{align}
\nonumber &\frac{\dd^k}{\dd\delta^k}\left|\nabla u(\delta)\right|^{h (\delta)}=\sum_{m=0}^{\infty} \frac{1}{m!} \sum_{\ell=0}^{k}\binom{k}{\ell}\\
&\times \left(\frac{\dd^{\ell}}{\dd\delta^{\ell}} h^{m}(\delta)\right)\frac{\dd^{k-\ell}}{\dd\delta^{k-\ell}}\left(\ln\left|\nabla u(\delta)\right|\right)^{m}\, .
\end{align}

We can now make a few observations.  First, given that $h\left(0\right)=0$, we must have that $h (\delta)\sim \mathcal{O}(\delta)$ as $\delta\to 0$.  Hence, $h^m(\delta)\sim \mathcal{O}\left(\delta^m\right)$ and so $\frac{\dd^{\ell}}{\dd\delta^{\ell}} h^{m}(\delta)\Big|_{\delta=0}\to 0 $ whenever $m>\ell$, which lets us restrict the upper range of $m$ and exchange the order of summation.  Further, when $m=0$ and $\ell\ne 0$, $\frac{\dd^{\ell}}{\dd\delta^{\ell}} h^{m}(\delta)$ clearly vanishes.  The one case in which $\ell=0$ and $m=0$ does not vanish is when $k=0$, which can be separated off on its own.  For all $k\ge 1$, we have
\begin{align}
\nonumber &\frac{\dd^k}{\dd\delta^k}\left|\nabla u(\delta)\right|^{h (\delta)}=\sum_{\ell=1}^{k}\binom{k}{\ell}\sum_{m=1}^{\ell} \frac{1}{m!} \left(\frac{\dd^{\ell}}{\dd\delta^{\ell}} h^{m}(\delta)\right)\\
&\times \frac{\dd^{k-\ell}}{\dd\delta^{k-\ell}}\left(\ln\left|\nabla u(\delta)\right|\right)^{m}\, .
\end{align}
To make further progress, we now leverage the Fa\`{a} di Bruno formula Eq.~\eqref{eq:diBrunos}, and find
\begin{align}
\nonumber \frac{1}{n!}\frac{\dd^n}{\dd\delta^n}\left(h (\delta)\right)^{m}&= \sum_{\substack{ \mathbf{p}\in S_n \\ \|\mathbf{p}\|_1 \le m}}\alpha_{m}(\mathbf{p})\left(h (\delta)\right)^{m-\|\mathbf{p}\|_{1}}\\
&\times \prod_{j=1}^{n}\left(\frac{h^{(j)}(\delta)}{j!}\right)^{p_j}\, ,
\end{align}
where
\begin{align}
\alpha_{m}(\mathbf{p})&=\frac{\prod_{p=0}^{\|\mathbf{p}\|_1-1}\left(m-p\right)}{\prod_{j=1}^{n}p_j!}\, .
\end{align}
Noting again the property that $h(0)=0$, we must have that $\|\mathbf{p}\|_1=m$, which yields
\begin{align}
\frac{\dd^n}{\dd\delta^n}\left(h (\delta)\right)^{m}&=n!\sum_{\substack{ \mathbf{p}\in S_n \\ \|\mathbf{p}\|_1 = m}}m!\prod_{j=1}^{n}\left(\frac{h^{(j)}(\delta)}{j!p_j!}\right)^{p_j}\, ,
\end{align}
Similarly, leveraging Fa\`{a} di Bruno we find
\begin{align}
\nonumber \frac{1}{n!}\frac{\dd^n}{\dd\delta^n}\left(\ln \left|\nabla u(\delta)\right|\right)^{m}&=\sum_{\substack{ \mathbf{p}\in S_n \\ \|\mathbf{p}\|_1 \le m}}\alpha_{m}(\mathbf{p})\left(\ln \left|\nabla u(\delta)\right|\right)^{m-\|\mathbf{p}\|_{1}}\\
&\times \prod_{j=1}^{n}\left(\frac{1}{j!}\frac{\dd^j}{\dd\delta^j}\ln \left|\nabla u(\delta)\right|\right)^{p_j}\, .
\end{align}
An additional application of Fa\`{a} di Bruno formula to the logarithm term gives
\begin{align}
\nonumber \frac{1}{j!}\frac{\dd^j}{\dd\delta^j}\ln \left|\nabla u(\delta)\right|\Big|_{\delta=0}&=\frac{1}{2}\sum_{\mathbf{r}\in S_j}\frac{\beta(\mathbf{r})}{\left|\nabla u_0\right|^{2\|\mathbf{r}\|_1}}\\
&\times \prod_{i=1}^{j}\left(\frac{1}{i!}\frac{\dd^i}{\dd\delta^i}\left|\nabla u(\delta)\right|^2\Big|_{\delta=0}\right)^{r_i}\, ,
\end{align}
where
\begin{align}
\beta(\mathbf{r})&=\frac{(-1)^{\|\mathbf{r}\|_1 -1}\left(\|\mathbf{r}\|_1-1\right)!}{\prod_{i=1}^{j}r_i!}\, .
\end{align}
Now, 
\begin{align}
\frac{\dd^i}{\dd\delta^i}\left|\nabla u(\delta)\right|^2\Big|_{\delta=0}&=\sum_{a=0}^{i}\binom{i}{a} \langle \frac{\dd^{i-a}}{\dd\delta^{i-a}}\nabla u(\delta),\frac{\dd^a}{\dd\delta^a}\nabla u(\delta)\rangle\Big|_{\delta=0}\, ,\\
&=i!\sum_{a=0}^{i}\langle \nabla u_{i-a},\nabla u_{a}\rangle\, ,
\end{align}
where the brackets here indicate an inner product over the dimensions of space (i.e., over the vector index of $\nabla$).  Hence,
\begin{widetext}
\begin{align}
\frac{1}{j!}\frac{\dd^j}{\dd\delta^j}\ln \left|\nabla u(\delta)\right|\Big|_{\delta=0}&=\frac{1}{2}\sum_{\mathbf{r}\in S_j}\frac{\beta(\mathbf{r})}{\left|\nabla u_0\right|^{2\|\mathbf{r}\|_1}}\prod_{i=1}^{j}\left(\sum_{a=0}^{i}\langle \nabla u_{i-a},\nabla u_{a}\rangle\right)^{r_i}\, ,\\
\Rightarrow \frac{1}{n!}\frac{\dd^n}{\dd\delta^n}\left(\ln \left|\nabla u(\delta)\right|\right)^{m}&=\sum_{\substack{ \mathbf{p}\in S_n \\ \|\mathbf{p}\|_1 \le m}}\alpha_{m}(\mathbf{p})\left(\ln \left|\nabla u_0\right|\right)^{m-\|\mathbf{p}\|_{1}}\prod_{j=1}^{n}\left(\frac{1}{2}\sum_{\mathbf{r}\in S_j}\frac{\beta(\mathbf{r})}{\left|\nabla u_0\right|^{2\|\mathbf{r}\|_1}}\prod_{q=1}^{j}\left(\sum_{a=0}^{q}\langle \nabla u_{q-a},\nabla u_{a}\rangle\right)^{r_q}\right)^{p_j}\, ,\\
%
\Rightarrow 
\nonumber \frac{\dd^k}{\dd\delta^k}\left|\nabla u(\delta)\right|^{h (\delta)}&=\sum_{\ell=1}^{k}k! \sum_{m=1}^{\ell} \sum_{\substack{ \mathbf{p}\in S_k \\ \|\mathbf{p}\|_1 = m}}\prod_{j=1}^{\ell}\left(\frac{h^{(j)}(\delta)}{j!p_j!}\right)^{p_j}\\
&\times \sum_{\substack{ \mathbf{p}\in S_{k-\ell} \\ \|\mathbf{p}\|_1 \le m}}\alpha_{m}(\mathbf{p})\left(\ln \left|\nabla u_0\right|\right)^{m-\|\mathbf{p}\|_{1}}\prod_{j=1}^{k-\ell}\left(\frac{1}{2}\sum_{\mathbf{r}\in S_j}\frac{\beta(\mathbf{r})}{\left|\nabla u_0\right|^{2\|\mathbf{r}\|_1}}\prod_{q=1}^{j}\left(\sum_{a=0}^{q}\langle \nabla u_{q-a},\nabla u_{a}\rangle\right)^{r_q}\right)^{p_j} \, .
\end{align}
\end{widetext}
Putting this all together, we find
\begin{widetext}
\begin{align}
\nonumber f(\delta)&=\sum_{n\ge 0}  \delta^n \nabla u_n+\sum_{n\ge 1}  \delta^n \sum_{k=1}^{n} \nabla u_{n-k} \sum_{\ell=1}^{k}  \sum_{m=1}^{\ell} \left(\sum_{\substack{ \mathbf{p}\in S_k \\ \|\mathbf{p}\|_1 = m}}\prod_{j=1}^{\ell}\left(\frac{h^{(j)}(\delta)}{j!p_j!}\right)^{p_j}\right)\\
&\times \sum_{\substack{ \mathbf{p}\in S_{k-\ell} \\ \|\mathbf{p}\|_1 \le m}}\alpha_{m}(\mathbf{p})\left(\ln \left|\nabla u_0\right|\right)^{m-\|\mathbf{p}\|_{1}}\prod_{j=1}^{k-\ell}\left(\frac{1}{2}\sum_{\mathbf{r}\in S_j}\frac{\beta(\mathbf{r})}{\left|\nabla u_0\right|^{2\|\mathbf{r}\|_1}}\prod_{q=1}^{j}\left(\sum_{a=0}^{q}\langle \nabla u_{q-a},\nabla u_{a}\rangle\right)^{r_q}\right)^{p_j} 
\end{align}
\end{widetext}

	\bibliographystyle{apsrev4-1}
	\bibliography{SeriesPaperRefs} 

@preamble{ " \newcommand{\noop}[1]{} " }

@article{HormanderALPDOII,
author={{H\"{o}rmander, L.}},
title={{The Analysis of linear partial differential operators II}},
journal={{Grundlehren}},
publisher={{Springer-Verlag}},
year={{1983}},
volume={{257}},
doi = {{https://doi.org/10.1007/b138375}},
}

@article{BIANCHINI2025128761,
title = {Existence and blow-up for non-autonomous scalar conservation laws with viscosity},
journal = {Journal of Mathematical Analysis and Applications},
volume = {542},
number = {1},
pages = {128761},
year = {2025},
issn = {0022-247X},
doi = {https://doi.org/10.1016/j.jmaa.2024.128761},
url = {https://www.sciencedirect.com/science/article/pii/S0022247X24006838},
author = {Stefano Bianchini and Giacomo Maria Leccese}
}

@article{GUIDOLIN2022126361,
title = {Global existence results for solutions of general conservative advection-diffusion equations in R},
journal = {Journal of Mathematical Analysis and Applications},
volume = {515},
number = {1},
pages = {126361},
year = {2022},
issn = {0022-247X},
doi = {https://doi.org/10.1016/j.jmaa.2022.126361},
url = {https://www.sciencedirect.com/science/article/pii/S0022247X22003754},
author = {P.L. Guidolin and L. Schütz and J.S. Ziebell and J.P. Zingano}
}

@article{KinnunenParviainen2010,
url = {https://doi.org/10.1515/acv.2010.002},
title = {Stability for degenerate parabolic equations},
author = {Juha Kinnunen and Mikko Parviainen},
pages = {29--48},
volume = {3},
number = {1},
journal = {Advances in Calculus of Variations},
doi = {doi:10.1515/acv.2010.002},
year = {2010},
lastchecked = {2025-10-22}
}

@article{fujita1966blowing,
  title={On the blowing up of solutions of the Cauchy problem for $u_t = \Delta u + u^{1 + \sigma}$},
  author={Fujita, Hiroshi},
  journal={J. Fac. Sci. Univ. Tokyo},
  volume={13},
  pages={109--124},
  year={1966}
}

@article{Bandle1994,
author = {Catherine Bandle and Howard A. Levine},
title = {{Fujita type phenomena for reaction-diffusion equations with convection like terms}},
volume = {7},
journal = {Differential and Integral Equations},
number = {5-6},
publisher = {Khayyam Publishing, Inc.},
pages = {1169 -- 1193},
year = {1994},
doi = {10.57262/die/1369329510},
URL = {https://doi.org/10.57262/die/1369329510}
}

@InProceedings{Friedman1988,
author="Friedman, Avner",
editor="Ni, W.-M.
and Peletier, L. A.
and Serrin, James",
title="Blow-up of solutions of nonlinear parabolic equations",
booktitle="Nonlinear Diffusion Equations and Their Equilibrium States I",
year="1988",
publisher="Springer New York",
address="New York, NY",
pages="301--318",
isbn="978-1-4613-9605-5"
}

@article{BANDLE19983,
title = {Blowup in diffusion equations: A survey},
journal = {Journal of Computational and Applied Mathematics},
volume = {97},
number = {1},
pages = {3-22},
year = {1998},
issn = {0377-0427},
doi = {https://doi.org/10.1016/S0377-0427(98)00100-9},
url = {https://www.sciencedirect.com/science/article/pii/S0377042798001009},
author = {Catherine Bandle and Hermann Brunner}
}

@article{barenblatt1952self,
  title={On self-similar motions of a compressible fluid in a porous medium},
  author={Barenblatt, Grigory I},
  journal={Akad. Nauk SSSR. Prikl. Mat. Meh},
  volume={16},
  number={6},
  pages={79--6},
  year={1952}
}

@BOOK{Coddington1956theory,
	title     = "Theory of ordinary differential equations",
	author    = "Coddington, Earl A and Levinson, Norman",
	publisher = "McGraw Hill Book Company",
	year      =  1955,
	address   = "New York"
}

@article{liao2009notes,
  title={Notes on the homotopy analysis method: some definitions and theorems},
  author={Liao, Shijun},
  journal={Communications in Nonlinear Science and Numerical Simulation},
  volume={14},
  number={4},
  pages={983--997},
  year={2009},
  publisher={Elsevier}
}

@article{LEE2006389,
title = {Large-time geometric properties of solutions of the evolution p-Laplacian equation},
journal = {Journal of Differential Equations},
volume = {229},
number = {2},
pages = {389-411},
year = {2006},
issn = {0022-0396},
doi = {https://doi.org/10.1016/j.jde.2005.07.028},
url = {https://www.sciencedirect.com/science/article/pii/S0022039606002452},
author = {Ki-ahm Lee and Arshak Petrosyan and Juan Luis Vázquez}
}

@article{LINDQVIST198793,
title = {Stability for the solutions of div(¦▽u¦p − 2▽u) = f with varying p},
journal = {Journal of Mathematical Analysis and Applications},
volume = {127},
number = {1},
pages = {93-102},
year = {1987},
issn = {0022-247X},
doi = {https://doi.org/10.1016/0022-247X(87)90142-9},
url = {https://www.sciencedirect.com/science/article/pii/0022247X87901429},
author = {Peter Lindqvist}
}

@article{dibenedetto1989cauchy,
  title={On the Cauchy problem and initial traces for a degenerate parabolic equation},
  author={DiBenedetto, Emmanuele and Herrero, Miguel A},
  journal={Transactions of the American Mathematical Society},
  volume={314},
  number={1},
  pages={187--224},
  year={1989}
}

@book{lindqvist2019notes,
  title={Notes on the stationary p-Laplace equation},
  author={Lindqvist, Peter},
  year={2019},
  publisher={Springer}
}

@article{jeng1969forced,
  title={Forced model equation for turbulence},
  author={Jeng, Dah-Teng},
  journal={The Physics of Fluids},
  volume={12},
  number={10},
  pages={2006--2010},
  year={1969},
  publisher={AIP Publishing}
}

@article{bec2007burgers,
  title={Burgers turbulence},
  author={Bec, J{\'e}r{\'e}mie and Khanin, Konstantin},
  journal={Physics reports},
  volume={447},
  number={1-2},
  pages={1--66},
  year={2007},
  publisher={Elsevier}
}

@article{okamura1983steady,
  title={Steady solutions of forced Burgers equation},
  author={Okamura, Makoto and Kawahara, Takuji},
  journal={Journal of the Physical Society of Japan},
  volume={52},
  number={11},
  pages={3800--3806},
  year={1983},
  publisher={The Physical Society of Japan}
}

@book{evans2022partial,
  title={Partial differential equations},
  author={Evans, Lawrence C},
  volume={19},
  year={2022},
  publisher={American mathematical society}
}

@article{aronsson1988p,
  title={On p-harmonic functions in the plane and their stream functions},
  author={Aronsson, Gunnar and Lindqvist, Peter},
  journal={Journal of Differential Equations},
  volume={74},
  number={1},
  pages={157--178},
  year={1988},
  publisher={Academic Press}
}

@article{dalzell2023quantum,
	author = {Dalzell, Alexander M and McArdle, Sam and Berta, Mario and Bienias, Przemyslaw and others},
	doi = {10.48550/arxiv.2310.03011},
	eprint = {2310.03011},
	journal = {arXiv},
	otherauthor = {Chen, Chi-Fang and Gily{\'e}n, Andr{\'a}s and Hann, Connor T and Kastoryano, Michael J and Khabiboulline, Emil T and Kubica, Aleksander and Salton, Grant and Wang, Samson and Brand{\~a}o, Fernando G S L},
	title = {{Quantum algorithms: A survey of applications and end-to-end complexities}},
	year = {2023}}

@article{morales2024quantum,
	author = {Morales, Mauro E S and Pira, Lirand{\"e} and Schleich, Philipp and Koor, Kelvin and others},
	doi = {10.48550/arxiv.2411.02522},
	eprint = {2411.02522},
	journal = {arXiv},
	otherauthor = {Costa, Pedro C S and An, Dong and Aspuru-Guzik, Al{\'a}n and Lin, Lin and Rebentrost, Patrick and Berry, Dominic W},
	title = {{Quantum Linear System Solvers: A Survey of Algorithms and Applications}},
	year = {2024}}

@article{harrowquantum2009,
	author = {Harrow, Aram W. and Hassidim, Avinatan and Lloyd, Seth},
	doi = {10.1103/physrevlett.103.150502},
	eprint = {0811.3171},
	issn = {0031-9007},
	journal = {Physical Review Letters},
	number = {15},
	pages = {150502},
	pmid = {19905613},
	title = {{Quantum Algorithm for Linear Systems of Equations}},
	volume = {103},
	year = {2009}}

@article{joseph2023quantum,
  title={Quantum computing for fusion energy science applications},
  author={Joseph, I and Shi, Y and Porter, MD and Castelli, AR and Geyko, VI and Graziani, FR and Libby, SB and DuBois, JL},
  journal={Physics of Plasmas},
  volume={30},
  number={1},
  year={2023},
  publisher={AIP Publishing}
}

@article{PhysRevResearch.2.043102,
  title = {Koopman--von Neumann approach to quantum simulation of nonlinear classical dynamics},
  author = {Joseph, Ilon},
  journal = {Phys. Rev. Res.},
  volume = {2},
  issue = {4},
  pages = {043102},
  numpages = {17},
  year = {2020},
  month = {Oct},
  publisher = {American Physical Society},
  doi = {10.1103/PhysRevResearch.2.043102},
  url = {https://link.aps.org/doi/10.1103/PhysRevResearch.2.043102}
}

@article{liu2021efficient,
	author = {Liu, Jin-Peng and Kolden, Herman {\O}ie and Krovi, Hari K. and Loureiro, Nuno F. and Trivisa, Konstantina and Childs, Andrew M.},
	doi = {10.1073/pnas.2026805118},
	eprint = {2011.03185},
	issn = {0027-8424},
	journal = {Proceedings of the National Academy of Sciences},
	number = {35},
	pages = {e2026805118},
	pmcid = {PMC8536387},
	pmid = {34446548},
	title = {{Efficient quantum algorithm for dissipative nonlinear differential equations}},
	volume = {118},
	year = {2021}}

@article{leyton2008quantum,
	author = {Leyton, Sarah K and Osborne, Tobias J},
	doi = {10.48550/arxiv.0812.4423},
	eprint = {0812.4423},
	journal = {arXiv},
	title = {{A quantum algorithm to solve nonlinear differential equations}},
	year = {2008}}

@article{krovi2023improved,
	author = {Krovi, Hari},
	doi = {10.48550/arxiv.2202.01054},
	eprint = {2202.01054},
	journal = {arXiv},
	title = {{Improved quantum algorithms for linear and nonlinear differential equations}},
	year = {2022}}

@article{bharadwaj2020quantum,
	author = {Bharadwaj, Sachin S and Sreenivasan, Katepalli R},
	doi = {10.29195/iascs.03.01.0015},
	journal = {Indian Academy of Sciences Conference Series},
	number = {1},
	title = {{Quantum Computation of Fluid Dynamics}},
	volume = {3},
	year = {2020}}

@article{jaksch-et-al-2023-variational-quantum-algorithms-for-computational-fluid-dynamics,
	author = {Jaksch, Dieter and Givi, Peyman and Daley, Andrew J and Rung, Thomas},
	doi = {10.2514/1.j062426},
	eprint = {2209.04915},
	issn = {0001-1452},
	journal = {AIAA Journal},
	number = {5},
	pages = {1885--1894},
	title = {{Variational Quantum Algorithms for Computational Fluid Dynamics}},
	volume = {61},
	year = {2023}}

@article{diening2013mini,
  title={Mini-Workshop: the p-Laplacian operator and applications},
  author={Diening, Lars and Lindqvist, Peter and Kawohl, Bernd},
  journal={Oberwolfach Reports},
  volume={10},
  number={1},
  pages={433--482},
  year={2013}
}

@incollection{frisch2002burgulence,
  title={Burgulence},
  author={Frisch, Uriel and Bec, J{\'e}r{\'e}mie},
  booktitle={New trends in turbulence Turbulence: nouveaux aspects: 31 July--1 September 2000},
  pages={341--383},
  year={2002},
  publisher={Springer}
}

@article{guermond2004mathematical,
  title={Mathematical perspectives on large eddy simulation models for turbulent flows},
  author={Guermond, J-L and Oden, J Tinsley and Prudhomme, Serge},
  journal={Journal of Mathematical Fluid Mechanics},
  volume={6},
  number={2},
  pages={194--248},
  year={2004},
  publisher={Springer}
}

@article{del2022finite,
  title={A finite difference method for the variational p-Laplacian},
  author={del Teso, F{\'e}lix and Lindgren, Erik},
  journal={Journal of Scientific Computing},
  volume={90},
  number={1},
  pages={67},
  year={2022},
  publisher={Springer}
}

@article{smagorinsky1963general,
  title={General circulation experiments with the primitive equations: I. The basic experiment},
  author={Smagorinsky, Joseph},
  journal={Monthly weather review},
  volume={91},
  number={3},
  pages={99--164},
  year={1963}
}

@article{hopf1950partial,
  title={The partial differential equation $u_t+uu_x=\mu_{xx}$},
  author={Hopf, Eberhard},
  journal={Communications on pure and applied mathematics},
  volume={3},
  number={3},
  pages={201--230},
  year={1950}
}

@article{cole1951quasi,
  title={On a quasi-linear parabolic equation occurring in aerodynamics},
  author={Cole, Julian D},
  journal={Quarterly of applied mathematics},
  volume={9},
  number={3},
  pages={225--236},
  year={1951}
}

@article{bender1991new,
  title={A new perturbative approach to nonlinear partial differential equations},
  author={Bender, Carl M and Boettcher, Stefan and Milton, Kimball A},
  journal={Journal of mathematical physics},
  volume={32},
  number={11},
  pages={3031--3038},
  year={1991},
  publisher={American Institute of Physics}
}
\end{document}